\documentclass[12pt]{extarticle}
\usepackage{amsmath, amsthm, amssymb, color}
\usepackage[colorlinks=true,linkcolor=blue,urlcolor=blue]{hyperref}
\usepackage{graphicx}
\usepackage{caption}
\usepackage{mathtools}
\usepackage{hyperref}
\usepackage{enumerate}
 \usepackage{comment}

\usepackage{cleveref}
\usepackage{verbatim}
\usepackage{tikz,tikz-cd,tikz-3dplot}
\usepackage{amssymb}
\usetikzlibrary{matrix}
\usetikzlibrary{arrows}
\usepackage{algorithm}
\usepackage[noend]{algpseudocode}
\usepackage{caption}
\usepackage[normalem]{ulem}
\usepackage{subcaption}
\usepackage{multicol}
\oddsidemargin \evensidemargin
\usepackage{makecell}
\usepackage{array}
\usepackage{enumitem}

 \newtheorem{theorem}{Theorem}[section]

\newtheorem{proposition}[theorem]{Proposition}
\newtheorem{lemma}[theorem]{Lemma}
\newtheorem{corollary}[theorem]{Corollary}

\newtheorem{conjecture}[theorem]{Conjecture}

\theoremstyle{definition}
\newtheorem{definition}[theorem]{Definition}
\newtheorem{problem}[theorem]{Problem}
\newtheorem{remark}[theorem]{Remark}

\newtheorem{example}[theorem]{Example}

\newcommand{\PP}{\mathbb{P}}
\newcommand{\RR}{\mathbb{R}}

\newcommand{\CC}{\mathbb{C}}
\newcommand{\ZZ}{\mathbb{Z}}

\title{\bf Graphical Scattering Equations}
\author{
  Barbara Betti,
  Viktoriia Borovik,
  Bella Finkel, \\
  Bernd Sturmfels and
  Bailee Zacovic}

\date{ }

\begin{document}
\maketitle

\begin{abstract} \noindent
The CHY scattering equations on the moduli space $\mathcal{M}_{0,n}$  play a prominent role at~the
interface of particle physics and algebraic statistics.  We study the scattering correspondence
when the Mandelstam invariants are restricted to a fixed graph on $n$ vertices.
\end{abstract}

\section{Introduction}

The {\em CHY scattering potential} in particle physics \cite{CHY} is the following function in $n$ unknowns:
\begin{equation}
\label{eq:scatpot} L(x) \,\,\, = \,\,\sum_{1\leq i < j \leq n} \!\! s_{ij} \,{\rm log}(x_i - x_j) . 
\end{equation}
The coefficients $s_{ij}$ are parameters, known as {\em Mandelstam invariants}.
They are assumed to satisfy $n$ {\em momentum conservation relations}.
Setting  $s_{ii} = 0$ and $s_{ij} = s_{ji}$, these relations are
\begin{equation}
\label{eq:momcon0}
\quad \sum_{j=1}^n s_{ij} \,=\, 0 \qquad \quad {\rm for} \,\, i = 1,2,\ldots,n.
 \end{equation}
 The critical points of the potential $L(x)$ are the solutions to the
 {\em CHY scattering equations}
 \begin{equation}
 \label{eq:scateqns} \qquad \frac{\partial L}{\partial x_i} \,\, = \,\,
 \sum_{j=1}^n \frac{s_{ij}}{x_i - x_j} \,\,\,= \,\,\,0 
 \qquad \quad {\rm for} \,\, i = 1,2,\ldots,n.  
 \end{equation}
 
The hypothesis (\ref{eq:momcon0}) ensures that (\ref{eq:scatpot})
is invariant under linear fractional transformations
\begin{equation}
\label{eq:PGL2b} x_i \, \,\mapsto \,\, \frac{a \,x_i + b}{c \,x_i +d } \qquad {\rm for}\, \,\, i = 1,2,\ldots,n.
\end{equation}
Here $x_i$ is a point in $\PP^1 = \CC \cup \{\infty\}$.
The solution set to (\ref{eq:scateqns}) consists of $(n-3)!$ orbits
under the ${\rm PGL}(2)$ action  in \eqref{eq:PGL2b}.
Each orbit is a point in the moduli space $\mathcal{M}_{0,n} \subseteq (\PP^1)^n/\!/{\rm PGL}(2)$.
 Hence, the scattering potential $L(x)$ in \eqref{eq:scatpot} has
$(n-3)!$ critical points in  $\mathcal{M}_{0,n}$;
see \cite{CHY, Lam, ST}.
Each of them corresponds to a surface of degree $n-2$ in $\PP^{n-1}$, given by an ideal as in (\ref{eq:I4}).

In this paper we ask  the following question: what happens when
some  $s_{ij}$ are set to zero?
Will the solutions to (\ref{eq:momcon0}) and (\ref{eq:scateqns}) still be reasonable?
This is the case in the study of
{\em minimal kinematics} in \cite{EPS}. We here offer a vast generalization.
Let us begin with two small scenarios.

\begin{example}[$n=4$] Assuming $x_i \not= x_j$, the equations (\ref{eq:momcon0}) and (\ref{eq:scateqns}) 
are equivalent to \vspace{-0.2cm}
$$ \begin{matrix}
& s_{12}-s_{34}\,\,=\,\, s_{13}-s_{24}\,\,=\,\, s_{14}-s_{23}\,\,=\,\, s_{12} + s_{13} + s_{14} & = & 0 \\
{\rm and} \quad &
s_{34} x_1 x_2+s_{24} x_1 x_3+s_{14}x_2x_3+s_{23} x_1 x_4+s_{13} x_2 x_4+s_{12} x_3 x_4 & = & \,0.
\end{matrix} \vspace{-0.2cm}
$$
The quadric defines a surface in $\PP^3$. It encodes
 a ${\rm PGL}(2)$ orbit in $(\PP^1)^4$, which is one point
 in $\mathcal{M}_{0,4}$.
If we set one parameter $s_{ij}$ to zero,
the system collapses, with no solutions in~$\mathcal{M}_{0,4}$.
\end{example}

\begin{example}[$n=5$]  \label{ex:n5three}
The scattering potential  has two critical points
on the surface $\mathcal{M}_{0,5}$.
There is one critical point if $s_{45} = 0$,
and also if $s_{23} = s_{45} = 0$.
On the other hand, the system collapses for $s_{34} = s_{45} = 0$.
There are no solutions 
when three $s_{ij}$ are zero. See Example~\ref{ex:ideals5}.
\end{example}
To address our question, we fix an undirected graph 
$G \subseteq \binom{[n]}{2}$ with no loops or multiple edges and with vertices in $[n] = \{1,\ldots,n\}$,
and we set $s_{ij} = 0$ for all non-edges $ij \not\in G$.
Then (\ref{eq:momcon0}) is a system of $n$ linear equations
in the variables $s_{ij}$ for $ij \in G$. Its solution space
in $\PP^{|G| -1}$ is the {\em kinematic space}, denoted $\mathcal{K}_G$.
The equations (\ref{eq:scateqns}) are linear in $s_{ij}$, and
they define a vector bundle over an open subset of~$\PP^{n-1}$. Its total space $\mathcal{V}_G$ is called 
the {\em scattering correspondence} of $G$. 
Its Zariski closure is an irreducible variety $\overline{\mathcal{V}}_G$ 
 in $\mathcal{K}_G \times \PP^{n-1}$. 
We want the graph $G$ to  have enough edges so that
the scattering correspondence $\mathcal{V}_G$ is non-empty.  Moreover, if we consider $\mathcal{V}_G$ modulo the  ${\rm PGL}(2)$ action \eqref{eq:PGL2b}, the two
natural projections
\begin{equation}
\label{eq:twoprojections} \mathcal{K}_G \,\leftarrow\, \mathcal{V}_G\, \rightarrow \,\mathcal{M}_{0,n}
\end{equation}
should have the desirable properties described in Definition~\ref{def:geomcop}.
For the complete graph $G = K_n$, we have
${{\rm codim}(\mathcal{V}_{K_n}) = 2n-3}$,
the left map is $(n-3)!$-to-one, and  fibers on the
right are linear subspaces of codimension $n-3$ in the kinematic space $\mathcal{K}_G$.
The graph $G$ is \emph{copious} if it has enough
edges for its projections to behave as in
Definition~\ref{def:geomcop}.
By Examples \ref{ex:n5three} and \ref{ex:ideals5}, only three  graphs with $n=5$ vertices are copious.
We introduce four different approaches to what copious means.
Our main theorem states that three of these~four agree.

\begin{theorem} \label{conj:main}
For every graph $G$, the following three conditions are equivalent:
\begin{enumerate}
\item $G$ is geometrically copious, \hfill (Definition \ref{def:geomcop})  \vspace{-0.23cm}
\item $G$ is matroidally copious, \hfill (Definition \ref{def:matrcop})  \vspace{-0.23cm}
\item $G$ is algebraically copious. \hfill (Definition \ref{def:algcop})  \vspace{-0.05cm}
\end{enumerate}
If $G$ has a universal vertex, then each of these conditions implies
\begin{enumerate}
\item[4.] $G$ is topologically copious.
\hfill (Definition~\ref{def:topcop}) \vspace{-0.05cm}
\end{enumerate}
\end{theorem}

Developing these different notions of copious is our main thread.
One important feature of the present paper is that
all $n$ variables $x_1,x_2,\ldots,x_n$ are treated equally.
We do not break their symmetry. In most other studies with explicit
 computations on $\mathcal{M}_{0,n}$, three of the points are fixed at
$0,1,\infty$. That choice, called \emph{gauge fixing}, is incompatible
with our graph $G$, unless $G$ has a {\em universal vertex}
(i.e., a vertex that is connected to all other vertices).
For this reason, our commutative algebra takes place
over a polynomial ring in $x_1,x_2,\ldots,x_n$.

We now discuss the organization of this article and its other main results.
In Section \ref{sec2} we examine the implications of momentum conservation.
Proposition \ref{prop: compatible_graph} rules out all graphs $G$
whose even cycle matroid has a co-loop. Lemma \ref{lem:isinvariant} leads to an 
algebraic formulation of our problem which 
respects the ${\rm PGL}(2)$ action on the $n$ variables in (\ref{eq:PGL2b}).
In Section \ref{sec3} we take a closer look at the scattering correspondence $\mathcal{V}_G$.
Its vector bundle structure defines the \emph{scattering matroid} $\mathcal{S}(G)$, which is 
our object of study in
Lemma   \ref{lem:scatmat} and Proposition \ref{prop:cand}.
The section concludes with the proof that geometrically copious equals matroidally copious.

In Section \ref{sec4} we turn to commutative algebra.
Building on recent advances in  \cite{KKMSW, Muh},
we develop tools for computing the prime ideal $I_G$ of the scattering correspondence.
This involves the equivalent concept of algebraically copious (Lemma \ref{lem: geom. iff alg. copious}).
We introduce the multidegree of the ideal $I_G$, and show in Proposition \ref{prop:acograph}
how this reveals
the \emph{maximum likelihood degree} $\mu(G)$.
Proposition \ref{prop:14} offers a census of 
copious graphs with $n=6$~vertices.

 In Section \ref{sec5} we define a very affine variety $\mathcal{M}_G$
 which is a partial compactification of the moduli space $\mathcal{M}_{0,n}$
 associated with the graph $G$.
This parallels the tropical construction due to Fry \cite{Fry2022}
on the moduli space of graphically stable rational curves. 
If $\mathcal{M}_G$ is smooth, then its signed Euler characteristic is the ML degree $\mu(G)$.
This rests on Huh's Theorem in~\cite{Huh} and is the content of  Theorem \ref{thm:pacom}.
Section \ref{sec5} concludes with a proof of our main theorem (Theorem \ref{conj:main}) and establishes a converse whenever $\mathcal{M}_G$ has a positive ML degree $\mu(G)$.

Section~\ref{sec6} is a case study for the complete bipartite graph $G=K_{s,t}$. We give a
determinantal description of the  variety $\mathcal M_{K_{s,t}}$ (Theorem~\ref{thm:bipartite}) and relate its ML degree $\mu(K_{s,t})$ to that of a 
conditional independence model (Theorem~\ref{thm:jose}). Section \ref{sec7} addresses the combinatorics of the ML degree.
When $\mathcal{M}_G$ is smooth,  Conjecture \ref{thm: MLdeg_admissible} suggests a general formula.
This is compared with the standard approach,
using bounded regions and chromatic polynomials, which 
applies for a special family of  graphs that allow 
gauge fixing (Theorem~\ref{cor:gauge}).

Section \ref{sec8} is our contribution to experimental mathematics.
We present exhaustive computations for graphs up to $n=9$ vertices, using both
symbolic and numerical methods.
Theorem \ref{thm:copious-graphs} summarizes the findings.
Our database of copious graphs is made available on {\tt Zenodo} \cite{BBBBV}.
We illustrate how the data can be used to gain insights and generate conjectures. In addition, we study graphs $G$ arising from hypertrees
\cite{CT, EPS}. We examine whether such graphs $G$ are copious, 
and we explore scattering equations in the setting of
  \cite[Section 5]{EPS}.  
  Our favorite hypertree is the bipyramid
  (Proposition \ref{prop:bipyramidTopologicallyCopious}, Conjecture \ref{conj:bipyramid} and Proposition~\ref{thm:bipyramid}).

In Section \ref{sec9}, we generalize the logarithmic discriminant
of \cite[Section 8]{KKT} from $K_n$ to arbitrary graphs $G$. Computational results on its degree
$\delta(G)$ are given in Proposition~\ref{prop:deltadata}.
Ample discussions
guide our readers to a plethora of fruitful directions for 
future research.

\section{Momentum Conservation}\label{sec2}

Given a graph $G \subseteq \binom{[n]}{2}$, its kinematic space
$\mathcal{K}_G$ is the linear subspace of $\CC^{|G|}$ defined by 
\begin{equation}
\label{eq:momcon}
 \sum_{j \in G_i} s_{ij} \,=\, 0 \qquad {\rm for} \,\, i = 1,2,\ldots,n.
 \end{equation}
 Here $G_i$ is the set of vertices of $G$ that are incident to the vertex  $i$.
 As in the Introduction, we usually identify the
 kinematic linear space $\mathcal{K}_G$ with its
 image in the projective space $\PP^{|G|-1}$.
 The $n$ linear forms in (\ref{eq:momcon}) determine a matroid
 $M(\mathcal{K}_G)$ of rank at most $n$ on the ground set~$G$.
Namely,
$M(\mathcal{K}_G)$ is the matroid of the vector configuration
 $\bigl\{ e_i +   e_j \in \RR^n : ij \in G \bigr\}$.
  This matroid is known as the {\em even cycle matroid} of the graph $G$.
 Zaslavsky's classical theory~\cite{Zas} refers to
$M(\mathcal{K}_G)$ as the {\em signed
graphic matroid} of the graph $G$ for negatively signed edges.

The following result is well-known at the interface of
graph theory and matroid theory.

\begin{lemma} \label{lem:circuits}
The circuits of the even cycle matroid $\,M(\mathcal{K}_G)$ are \vspace{-0.2cm}
\begin{enumerate}
    \item even cycles; \vspace{-0.2cm}
    \item pairs of odd cycles intersecting in a single vertex; \vspace{-0.2cm}
    \item pairs of vertex-disjoint odd cycles connected by a simple path.
\end{enumerate}
\end{lemma}

We say that $G$ is {\em compatible with momentum conservation} if
$\mathcal{K}_G \cap (\CC^*)^{|G|}$ is non-empty.
This means that  momentum conservation 
can be achieved with $s_{ij}$ that are supported on $G$.

\begin{proposition}\label{prop: compatible_graph}
A graph $G$ is compatible with momentum conservation if and only if its even cycle matroid
$M(\mathcal{K}_G)$ has no co-loop.
This holds if and only if for every edge $e$ of $G$, the graph
$G \setminus e$ has the same number of bipartite connected components as $G$. 
\end{proposition}

\begin{proof}
The graph $G$ is not compatible with momentum conservation if and only if
$\mathcal{K}_G$ lies in some coordinate hyperplane $\{s_{ij} = 0\}$ of $\CC^{|G|}$.
This means that $ij$ is a co-loop of $M(\mathcal{K}_G)$.
The second statement can be derived from Lemma~\ref{lem:circuits}.
In Zaslavsky's framework, it follows from the fact that the
balanced subgraphs are  the bipartite subgraphs 
 \cite[Corollary~7.D.1]{Zas}.
 \end{proof}

Moreover, as a corollary of~\cite[Corollary~7.D.3(j)]{Zas}, we obtain the following statement.
\begin{corollary} \label{cor: bipartite connected components}
The codimension of $\mathcal{K}_G$ equals the rank of  $M(\mathcal{K}_G)$. This equals
$n$ minus~the number of bipartite connected components of $G$, where isolated vertices are 
viewed as bipartite.
\end{corollary}

\noindent We shall use the following abbreviation for the codimension of the kinematic space:
\begin{equation}
\label{eq:kappa}
 \kappa(G) \,\, \coloneqq \,\, {\rm codim}(\mathcal{K}_G) \,\, = \,\, {\rm rank}( M(\mathcal{K}_G)). 
 \end{equation}

Our object of study in  this paper is the \emph{graphical scattering potential}
\begin{equation*}\label{eq:scatpotential}
     L_G(x) \,\, = 
\,\,\sum_{ij \in G} s_{ij} \,{\rm log}(x_i - x_j).
\end{equation*}
We assume that the Mandelstam parameters $s_{ij}$, for $ij\in G$,
are nonzero and satisfy~\eqref{eq:momcon}. For such parameters to exist,
the graph $G$ must be compatible with momentum conservation.
\begin{lemma} \label{lem:isinvariant}
The scattering potential $ L_G$  is constant on the orbits of the
${\rm PGL}(2)$~action~\eqref{eq:PGL2b}.
Its differential ${\rm d}L_G$ is a well-defined
$\operatorname{PGL}(2)$-invariant logarithmic one-form on~$\mathcal{M}_{0,n}$.
\end{lemma}

\begin{proof}
We consider the lattice $  \mathcal{K}_G \,\cap \, \ZZ^{|G|}$
of all integer vectors in the kinematic space. This is a 
free abelian group of rank $r \coloneqq |G|-\kappa(G)$. 
Every vector $u \in  \mathcal{K}_G \,\cap \, \ZZ^{|G|}$
gives rise to a character $\tilde u : (\CC^*)^{|G|} \rightarrow \CC^*$.
We write $\tilde u$ as a Laurent monomial  in the
quantities $x_i - x_j$ for $ij \in G$. Then $\tilde u$
is a rational function on $\PP^{n-1}$ which is constant
on ${\rm PGL}(2)$ orbits.
 We call $\tilde u$ the {\em generalized cross-ratio},
 and we view it as a regular function on
 the moduli space $\mathcal{M}_{0,n}$.

Let $u_1,u_2,\ldots,u_r$ be a basis of 
 $  \mathcal{K}_G \,\cap \, \ZZ^{|G|}$. 
Since $s = (s_{ij})$ satisfies momentum conservation,
we can write
$s = t_1 u_1 + t_2 u_2 + \cdots + t_r u_r$, where
the $t_k$ are coordinates on $\mathcal{K}_G$.
We then have
\begin{equation}
\label{eq:LGinvariant}
 L_G(x) \,\, = \,\,\sum_{k=1}^r t_k \cdot {\rm log} (\tilde u_k(x)) . 
 \end{equation}
 This is an expression in terms of generalized cross-ratios, so it is
 ${\rm PGL}(2)$ invariant. As a consequence of~\eqref{eq:LGinvariant}, we also obtain
${\rm d}L_G=\sum_{k=1}^r t_k\,\frac{d\widetilde{u}_k}{\widetilde{u}_k}$,
which is well-defined on
$\mathcal{M}_{0,n}$. \end{proof}

The following graph appears prominently in
\cite[Section 5]{EPS} and in \cite[Theorem 5.2]{KKMSW}.

\begin{example}[Octahedron] \label{ex:octah}
Let $G$ be the graph with $n=6$ and non-edges $12,34,56$.
Here $r = 6$, and we choose the following integer basis for the kinematic subspace
$\mathcal{K}_G $ in $\CC^{|G|} = \CC^{12}$:
\begin{equation}
\label{eq:octabasis1}
 \begin{matrix}
u_1 \,=\,  e_{13} + e_{24} - e_{14} - e_{23}, &
u_2 \,=\,  e_{15} + e_{26} - e_{16} - e_{25},   &                                                  
u_3 \,=\, e_{35} + e_{46} - e_{36} - e_{45}, \\
u_4 \,=\, e_{14} + e_{25} - e_{15} - e_{24}, &
u_5 \,=\, e_{23} + e_{45} - e_{24} - e_{35}, &
u_6 \,=\, e_{25} + e_{46} - e_{26} - e_{45}.
\end{matrix}
\end{equation}
We can thus write  $s = \sum_{k=1}^6 t_k u_k$ for $s \in \mathcal{K}_G$.
The corresponding  cross-ratios~are
\begin{equation}
\label{eq:octabasis2}
 \begin{matrix}
\tilde u_1(x) \,=\,  \frac{(x_1-x_3)(x_2-x_4)}{(x_1-x_4)(x_2-x_3)},  &
\tilde u_2(x) \,=\, \frac{(x_1-x_5)(x_2-x_6)}{(x_1-x_6)(x_2-x_5)}, &
\tilde u_3(x) \,=\, \frac{(x_3-x_5)(x_4-x_6)}{(x_3-x_6)(x_4-x_5)}, \medskip \\
\tilde u_4(x) \,=\, \frac{(x_1-x_4)(x_2-x_5)}{(x_1-x_5)(x_2-x_4)}, &
\tilde u_5(x) \,=\, \frac{(x_2-x_3)(x_4-x_5)}{(x_2-x_4)(x_3-x_5)}, &
\tilde u_6(x) \,=\, \frac{(x_2-x_5)(x_4-x_6)}{(x_2-x_6)(x_4-x_5)}.
\end{matrix}
\end{equation}
In terms of these six cross-ratios,  the scattering potential on the octahedron graph $G$ equals
\begin{equation}
\label{eq:octaLG}
 L_G(x) \, = \, \sum_{ij \in G} s_{ij} \,{\rm log}(x_i-x_j) \, = \, \sum_{k=1}^6 t_k \,{\rm log}(\tilde u_k(x)) .
\end{equation}
This demonstrates how $L_G(x)$ can be written manifestly as a function
on $\mathcal{M}_{0,6}$, but without breaking the symmetry among the $n$ points.
As was said before, we abstain from gauge~fixing.
\end{example}

For any graph $G$ on $n$ vertices, we wish to solve  the
{\em graphical scattering equations}
 \begin{equation}
\label{eq:gradient}  \frac{\partial{L_G}}{\partial x_1} \,= \,
\frac{\partial{L_G}}{\partial x_2} \,= \,\,\cdots \,\,= \,
\frac{\partial{L_G}}{\partial x_n} \,= \, 0. 
\end{equation}
The {\em braid arrangement} $\mathcal{H}$ consists of the hyperplanes $\{x_i = x_j \}$ in $ \CC^{n}$ for $ij \in {[n] \choose 2}$. 
Let $u \in \CC^{n} \backslash \mathcal{H}$ be a solution of (\ref{eq:gradient}). Then
all points in the ${\rm PGL}(2)$ orbit of $u$ are
 solutions. This orbit is a threefold in $\CC^n$. Using $a,b,c,d$ 
 as unknowns,  the orbit has the parametrization
 \begin{equation}
\label{eq:PGL2} x_i \, = \, \frac{a u_i + b}{c u_i +d }, \qquad {\rm with}\, \,\, cu_i+d \neq 0 \qquad {\rm for}\, \,\, i = 1,2,\ldots,n. 
\end{equation}
Rewriting equations~\eqref{eq:PGL2} as
$c u_i x_i+d x_i-a u_i-b=0$ for $i= 1,\ldots,n$,
we see that $x$ lies in the closure of the
${\rm PGL}(2)$ orbit of $u$ precisely when the following matrix has rank at most~three:
$$
U \coloneqq \begin{pmatrix}
1&1&\cdots&1\\
u_1&u_2&\cdots&u_n\\
x_1&x_2&\cdots&x_n\\
u_1x_1&u_2x_2&\cdots&u_nx_n
\end{pmatrix}.
$$
 Thus, eliminating $a,b,c,d$
from~\eqref{eq:PGL2} yields the ideal $I_4(U)$ of its $4\times4$ minors. Equivalently, this 
ideal is
generated by the $2\times2$ minors below written in  the physics-inspired~notation:
\begin{equation} \label{eq:I4}
I_4 (U) \,= \,\,
I_2\! \begin{small} \begin{pmatrix}
[13]\langle 23 \rangle &  [14] \langle 24 \rangle & \cdots & [1n]\langle 2n\rangle \smallskip \\
[23]\langle 13 \rangle &  [24] \langle 14 \rangle & \cdots & [2n]\langle 1n\rangle \\
\end{pmatrix} \end{small}
\,\subset \,\CC[x], \quad \text{where}  \quad \begin{matrix} [ij] = x_i - x_j,
\\  \langle i j \rangle = u_i - u_j.
\end{matrix} 
\end{equation}

We now explain the second determinantal presentation. Let $U_i$ denote the $i$-th column of $U$. Since $u_1\neq u_2$, the columns $U_1$ and $U_2$ are linearly independent. Consider the linear~map
$$
\pi: \CC^4 \to \CC^2, \quad (1,u,x,ux)^{\top}\mapsto
((x_1-x)(u_2-u), \, 
(x_2-x)(u_1-u))^\top.
$$
Its kernel is $\operatorname{span}(U_1,U_2)$ and
$\pi(U_i)= ([1i]\langle2i\rangle,\,  [2i]\langle1i\rangle)^\top$  for each $i\geq3$.
Hence $U$ has rank at most three if and only if the images of
$U_3,\ldots,U_n$ in $\CC^2$
span a vector space of dimension at most one. Equivalently, all $2\times2$ minors of the $2 \times (n-2)$ matrix in~\eqref{eq:I4} vanish.

One can apply~\cite[Theorem~2.1]{EisenbudSections} to the matrix in~\eqref{eq:I4} to conclude that its ideal of $2\times2$ minors is prime.  That ideal is our algebraic encoding for
one  point $u$ in the moduli space $ \mathcal{M}_{0,n}$.
The variety defined by (\ref{eq:I4}) in $\PP^{n-1}$ is   a surface of degree $n-2$.
This surface is the projectivized affine chart of the orbit closure of the
point $(u_1,\ldots,u_n) \in (\PP^1)^n$.


We expect the scattering equations (\ref{eq:gradient}), for generic $s \in \mathcal{K}_G$,
to have a finite, non-empty set of solutions in~$\mathcal{M}_{0,n}$. 
This expectation will be made
precise in the next section.
The~number of solutions is the {\em ML degree} $\mu(G)$.
Thus, the variety defined by (\ref{eq:gradient}) in $\PP^{n-1}$~consists~of $\mu(G)$ surfaces of degree $n-2$, each given
by a  prime ideal (\ref{eq:I4}). For instance,  
for the~octahedron  $G$ from Example~\ref{ex:octah}, we find
 exactly $\mu(G) = 2$ such surfaces in $\PP^5$.
See {\cite[Example~5.4]{EPS}}.

We close this section by highlighting the
similarities and differences between  \cite{KKMSW} and the present paper.
 Let $\mathcal{H}_G$ be the \emph{graphical arrangement}, obtained from $\mathcal{H}$ by deleting all hyperplanes 
 with label $ij \notin G$.  
The likelihood correspondence of the arrangement 
$\mathcal{H}_G$ is studied in
\cite[Section 4]{KKMSW}. 
The similarity is that
both theories concern the same
function~$L_G$. 

The crucial difference lies in the assumptions on
the parameters $s_{ij}$. In \cite{KKMSW} one has
$\sum_{1 \leq i < j \leq n} s_{ij} = 0$, but otherwise the
$s_{ij}$ are unconstrained. This results in finitely
many critical points in $\PP^{n-1}$. By contrast, we here
assume momentum conservation (\ref{eq:momcon}).
This assumption implies that
 the critical points in $\PP^{n-1}$ are not
isolated, but they come in ${\rm PGL}(2)$ orbits.
These are the surfaces defined by  (\ref{eq:I4}), representing
points in the moduli space $\mathcal{M}_{0,n}$.

The likelihood ideal
in \cite{KKMSW} is usually
generated by M\"uhlherr's equations from \cite{Muh}.
We refer to \cite[Theorem 5.2]{KKMSW}, which says
that the octahedron is the only exception for~$n=6$.
M\"uhlherr's theorem in \cite{Muh} is  relevant here, and will be discussed in Section~\ref{sec4}.
However, the commutative algebra of our scattering ideals 
is certainly more complicated than that in \cite{KKMSW}. 

\section{The Scattering Matroid}\label{sec3}

In this section, we study the \emph{scattering correspondence} 
$\mathcal{V}_G$  given by (\ref{eq:momcon}) and (\ref{eq:gradient}). Namely, $\mathcal{V}_G$
   is the  set of pairs $(s,x)$ in $\PP^{|G|-1}\times \PP^{n-1}$
where $s$ ranges over generic points in $\mathcal{K}_G$ and
$x \in \PP^{n-1}\backslash \mathcal{H}$ is a critical point of
the graphical scattering potential $L_G(x)$. In symbols:
\begin{equation}
\label{eq:V_G}
\mathcal{V}_G \,\,=\,\,  \left\{ (s,x)\in \PP^{|G|-1}\times \PP^{n-1}\backslash \mathcal{H} :\, \frac{\partial L_G}{\partial x_i}=0
\,\,\,{\rm and} \,\, \sum\limits_{j\in G_i} s_{ij}=0 \ \text{ for }  \ i=1,\dots,n \right\}. 
\end{equation}
We write $\overline{\mathcal{V}}_G$  for the Zariski closure of $\mathcal{V}_G$ in
  $\mathcal{K}_G\times \PP^{n-1}$.
We frequently pass from $\PP^{n-1}\backslash \mathcal{H}$ to $\mathcal{M}_{0,n}$ by taking the quotient modulo ${\rm PGL}(2)$.
We thus view $\mathcal{V}_G$ as a subvariety of $\mathcal{K}_G \times \mathcal{M}_{0,n}$.

The scattering correspondence $\mathcal{V}_G$ 
is defined by rational function equations (\ref{eq:gradient}).
It is far from obvious how 
 these can be replaced with polynomials. How can we
clear denominators gracefully? This will be addressed in Section \ref{sec4}.
For a first glimpse of that issue,  consider
the following polynomials of bidegree $(1,k)$ in $(s,x)$.
They are symmetric functions in $(x_i,x_j)$:
\begin{equation} \label{eq:alwaysthere}
 \qquad \sum\limits_{ij\in G}\left( s_{ij} \sum\limits_{\ell=0}^{k} x_i^{k-\ell}x_j^\ell \right),
 \qquad {\rm where} \,\, k = 2,3,\ldots,n-2.
\end{equation} 
For $k=1$, this is a linear form in $x_1,\ldots,x_n$ whose 
coefficients are precisely the $n$ momentum conservation relations~\eqref{eq:momcon}.
Hence, it vanishes on $\mathcal{V}_G$. The same is true for $k \geq 2$ as well.

\begin{lemma} The polynomials in  (\ref{eq:alwaysthere}) vanish 
on the scattering correspondence $\mathcal{V}_G$.
\end{lemma}
\begin{proof}
Up to an index shift, these are the ideal generators $\theta_k$ in \cite[Theorem 5.8]{KKMSW}.
\end{proof}

We now come to the first definition of copious.
It refers to a graph with enough edges.
For many graphs, but not all, this property is preserved under adding edges. See
Proposition~\ref{prop: add edges}.

\begin{definition} \label{def:geomcop}
    A graph $G$ is \emph{geometrically copious} if the projection of 
    $\mathcal{V}_G$ to $\mathcal{M}_{0,n}$ is dominant and the projection to $\mathcal K_G$
is dominant with generic fiber of dimension two. In the formulation of (\ref{eq:twoprojections}), the latter condition means
    that the map on the left is finite-to-one.
\end{definition}
\noindent
If $G$ is geometrically copious then its ML degree $\mu(G)$ is
a well-defined positive integer. It is the degree of the
map on the left in (\ref{eq:twoprojections}),
i.e.,~the number of critical points of $L_G(x)$ in $\mathcal{M}_{0,n}$.

\begin{proposition}\label{prop:degree2_inadmissible}
    Every vertex in a geometrically copious graph $G$ has degree at least three.
\end{proposition}

\begin{proof}
Consider a vertex $i$ of degree $\leq 2$ in $G$.
    If $G_i=\{j,\ell\}$  then (\ref{eq:momcon}) implies $s_{ij}=-s_{i\ell}$, and 
     $\,{\partial L_G}/{\partial x_i} = 0\,$ forces $x_j=x_\ell$. Hence the map
     $\mathcal{V}_G\rightarrow \mathcal{M}_{0,n}$ is not dominant.
     If $G_i \subseteq \{\ell\}$ then the variable $x_i$ does not occur in the scattering equations, so the fibers of $\mathcal{K}_G \leftarrow \mathcal{V}_G$ are either empty or positive-dimensional.
    In both cases, $G$ is not geometrically copious.
\end{proof}

A key role in what makes a graph $G$ copious 
is played by matroids on the ground set $G$.
We  saw a first instance of this in Proposition \ref{prop: compatible_graph}.
These matroids furnish the combinatorial framework for $\mathcal{V}_G$
 as a vector bundle over  $\mathcal{M}_{0,n}$.
The point is that the equations (\ref{eq:momcon}) are linear over $\ZZ$,
and the equations (\ref{eq:gradient}) are linear over the field~$\CC(x)$.
We define the {\em scattering matroid} $\mathcal{S}(G)$ to be the matroid on $G$ that is
 given by the $2n$ linear equations in (\ref{eq:momcon}) and~(\ref{eq:gradient}).

\begin{lemma} \label{lem:scatmat}
The scattering matroid $\mathcal{S}(G)$ has rank at most $2n-3$.
This upper bound is attained for the complete graph $K_n$.
In symbols, we have ${\rm rank}(\mathcal{S}(K_n)) = 2n-3$.
\end{lemma}

\begin{proof}
The fiber in $\mathcal{V}_G$ over 
a generic point $u$ in $\PP^{n-1}$
is the projectivization of a subspace of the vector space $\CC^{|G|}$.
By definition, $\mathcal{S}(G)$ is the matroid defined by this subspace. 
We represent $\mathcal{S}(G)$ by the
following $2n \times |G|$ matrix $S_G$ with entries in the field $\CC(x_1,\dots, x_n)$.
The columns of  $S_G$ are indexed by edges $ij \in G$.
The rows of $S_G$ are the coefficients
 of (\ref{eq:momcon}) and (\ref{eq:gradient}) with respect to the $s_{ij}$.
 Explicitly, the entries of $S_G$ in a given column $ij $ are as follows.
 For $k=1,\ldots,n$,
 the entry in row $k$ is
 $1$ if $k \in \{i,j\}$ and it is $0$ otherwise.
 The entry in row $n+k$ is
 $(x_i-x_j)^{-1}$ if $i=k$, the entry  is
 $(x_j-x_i)^{-1}$ if $j=k$, and it is $0$ otherwise.
 
 Fix the complete graph $G = K_n$,
 and consider the $2n \times \binom{n}{2}$
 matrix $S_{K_n}$ defined above.
  Since $I_{K_n}$ is a complete intersection by \cite[Example 2.1]{KKMSW},
generated by (\ref{eq:momcon}) and the
$n-3$ polynomials in (\ref{eq:alwaysthere}),
we can  replace the last $n$ rows
of our matrix $S_{K_n}$ by the $n-3$ coefficient vectors
of (\ref{eq:alwaysthere}).
The resulting $(2n-3) \times \binom{n}{2}$ matrix
has linearly independent rows, and it
has the same row space as $S_{K_n}$.
It  therefore defines the  same matroid
$\mathcal{S}(K_n)$ of rank $2n-3$.

The matrix $S_G$ for any subgraph $G$ of $K_n$ is obtained from the
$(2n-3) \times \binom{n}{2}$ matrix $S_{K_n}$ by deleting all
columns corresponding to non-edges of $G$.
The rank of a matrix cannot increase if we delete
some columns. Hence $S_G$ has rank at most $2n-3$, and so does $\mathcal{S}(G)$.
\end{proof}

The proof  suggests that we start by examining the matrix $S_{K_n}$ for the
complete graph.
For any graph $G$, the matrix $S_G$ is obtained
 by deleting the columns of $S_{K_n}$ indexed by non-edges of $G$.
The codimension $\kappa(G)$ in \eqref{eq:kappa} is  the rank of the submatrix of $S_G$  given by the first $n$ rows.
For $G$ to be copious, the whole matrix $S_G$ should have rank  $\kappa(G)+n-3$.
\begin{example}[$n=6$]
The matrix $S_{K_6}$ has $12$ rows and $15$ columns, and it equals
$$ \!\!\!
\begin{footnotesize}
\begin{bmatrix}
 1 & 1 & 1 & 1 & 1 & 0 &  \cdots  & 0 & 0 \\
 1 & 0 & 0 & 0 & 0 & 1 & \cdots  & 0 & 0 \\
 0 & 1 & 0 & 0 & 0 & 1 & \cdots & 0 & 0 \\
 0 & 0 & 1 & 0 & 0 & 0 & \cdots & 1 & 0 \\
 0 & 0 & 0 & 1 & 0 & 0 &  \cdots & 0 & 1 \\
 0 & 0 & 0 & 0 & 1 & 0 &  \cdots & 1 & 1 \\ 
  (x_1{-}x_2)^{{-}1}\! & \! (x_1{-}x_3)^{{-}1}\! & \! \!(x_1{-}x_4)^{{-}1} \!&\!  (x_1{-}x_5)^{{-}1}\! & 
  \! \!(x_1{-}x_6)^{{-}1} \!& 
   0 & \cdots & 0 & 0 \\
  (x_2{-}x_1)^{{-}1} \!& 0 & 0 & 0 & 0 & \!\!\!(x_2{-}x_3)^{{-}1} &  \cdots & 0 & 0 \\
  0 &\! (x_3{-}x_1)^{{-}1}\! & 0 & 0 & 0 &\! \!\!(x_3{-}x_2)^{{-}1} & \cdots & 0 & 0 \\
  0 & 0 &\! \!(x_4{-}x_1)^{{-}1} \!& 0 & 0 & 0 & \cdots & (x_4{-}x_6)^{{-}1} & 0 \\
  0 & 0 & 0 &\! (x_5{-}x_1)^{{-}1}\! & 0 & 0 & \cdots & 0 &\!\!\! \!(x_5{-}x_6)^{{-}1} \\
  0 & 0 & 0 & 0 &\!\! (x_6{-}x_1)^{{-}1} \! & 0 & \cdots & (x_6{-}x_4)^{{-}1} & \!\!\!\! (x_6{-}x_5)^{-1} \\
 \end{bmatrix}.
 \end{footnotesize}
$$
\end{example}
\begin{remark}
The scattering matroid $\mathcal{S}(G)$ is Kalai's $2$-hyperconnectivity matroid~${\cal H}_2(G)$, see~\cite{Kal85}. Indeed, multiply the $ij$-column of $\mathcal{S}(G)$ by $x_i-x_j$ and perform the row~operations:
$$
L(i)\leftarrow x_iL(n+i)-L(i),\quad \text{for } \; i=1,\ldots,n.
$$
After grouping the rows by vertices, the $ij$-column has 
$p_j=(x_j,1)^\top$ at $i$ and $-p_i=-(x_i,1)^\top$ at  $j$,
which, up to transposition, is the hyperconnectivity matrix of the points
$p_i=(x_i,1)$. Thus ${\cal S}(G)={\cal H}_2(G)$. In particular, the
rank statement of Lemma~\ref{lem:scatmat} is the known identity
$\operatorname{rank}{\cal H}_2(K_n)=2n-3$.
We thank the anonymous referee for pointing out this connection.
\end{remark}
\begin{proposition} \label{prop:cand}
If $G$ is geometrically copious then
 the following constraints~hold:
 \begin{equation}
\label{eq:candidate} 
{\rm rank}\bigl( \mathcal{S}(G))  \,\,=\,\, {\rm rank}\bigl( \mathcal{S}(G \backslash ij )\bigr) \,\,=\,\, \kappa(G) \,+\,n-3
\qquad \hbox{for all edges} \,\, ij \in G,
\end{equation}  
\begin{equation}
\label{eq:fromxtou2} {\rm and} \qquad
\sum_{\substack{x\in\PP^{n-1}\setminus\mathcal H\\ x\ \mathrm{generic}}} \! {\rm ker}(S_G(x)) \,\, = \,\, \mathcal{K}_G. 
\end{equation}
The last condition is equivalent to
\begin{equation}
\label{eq:fromxtou1}
 {\rm dim} \,\bigl(\bigcap_{\substack{x\in\PP^{n-1}\setminus\mathcal H\\ x\ \mathrm{generic}}} {\rm rowspace}(S_G(x))\, \bigr) \,\, = \,\, \kappa(G).
\end{equation}
\end{proposition}

\begin{proof}
    We work in the setting of (\ref{eq:twoprojections}), where
$ \mathcal{V}_G \subset \mathcal{K}_G \times \mathcal{M}_{0,n}$.
Our hypothesis says that
the map to $\mathcal{K}_G$ is dominant and finite-to-one. This means that the fiber over a generic point in $\mathcal{M}_{0,n}$
does not lie in a proper subspace of $\mathcal{K}_G$. In particular, that fiber doesn't lie in a coordinate subspace $\{s_{ij} = 0\}$, that is, no edge $ij$ of $G$ is a co-loop of the scattering matroid $\mathcal{S}(G)$. In symbols, this is the first equality of (\ref{eq:candidate}).
The fact that the map $\mathcal{V}_G \to \mathcal{K}_G$ is of finite degree also implies that the codimension of the
scattering correspondence in $\mathcal{K}_G \times \mathcal{M}_{0,n}$
equals $n-3 =  {\rm dim}(\mathcal{M}_{0,n})$.
Its codimension in $\PP^{|G|-1} \times \mathcal{M}_{0,n}$ is therefore
${\rm codim}(\mathcal{K}_G) +n-3 = \kappa(G) + n-3$.
The latter codimension is precisely the rank of 
$\mathcal{S}(G)$, thus  \eqref{eq:candidate} holds.

For any $x \in \PP^{n-1}\backslash \mathcal{H}$,
the kernel of $S_G(x)$ is contained in $\mathcal{K}_G$, because $\mathcal{K}_G$
 is the kernel of the first $n$ rows of $S_G(x)$. Hence
the left hand side in (\ref{eq:fromxtou2}) is a subspace of $\mathcal{K}_G$. The fact that the map to $\mathcal{K}_G$ is dominant implies that 
there are no points of $\mathcal{K}_G$ outside that subspace.
Hence, the equivalent conditions \eqref{eq:fromxtou2}–\eqref{eq:fromxtou1} hold for any geometrically copious graph $G$.
\end{proof}
Our second notion of what it means for a graph $G$ to be copious rests on matroid theory.

\begin{definition} \label{def:matrcop}
A graph $G\!\subseteq\!\binom{[n]}{2}$ is  called \emph{matroidally copious} 
if both (\ref{eq:candidate}) and (\ref{eq:fromxtou2})$\simeq$(\ref{eq:fromxtou1}) are~satisfied. In particular,    the rank of the scattering matroid $\mathcal{S}(G)$ satisfies the upper bound
    \begin{equation}\label{eq:ScatMatrank}
        \kappa(G) + n-3 < |G|.
    \end{equation} 
        \end{definition}      
\begin{theorem}\label{lem:matroidally_geometrically}
    A graph is matroidally copious if and only if it is geometrically copious.
\end{theorem}

\begin{proof}
Proposition~\ref{prop:cand} proves the `if' direction. Conversely, suppose
that $G$ is matroidally copious. We first show that the map
$\mathcal V_G\rightarrow\mathcal K_G$ is dominant. As in Lemma~\ref{lem:isinvariant}, we set $r=|G|-\kappa(G) = \dim\mathcal K_G$, and we choose a $\ZZ$-basis $u_1,\ldots,u_r$ of the kinematic~lattice~${\mathcal K_G\cap\ZZ^{|G|}}$. The corresponding generalized
cross-ratios $z_1=\widetilde u_1(x), \ldots, z_r=\widetilde u_r(x)$
define the~map 
$$\varphi:\mathcal M_{0,n}\rightarrow(\CC^*)^r, \quad x\mapsto
\bigl(\widetilde u_1(x),\ldots,\widetilde u_r(x)\bigr).$$
Let $\mathcal M_G = \overline{{\rm im}(\varphi)}$ and, for $s=\sum_{i=1}^r t_i u_i\in\mathcal K_G$, take the invariant
logarithmic one-form 
\begin{equation}\label{eq: 1-form}
    \omega_s\,=\,\sum_{i=1}^r t_i\,\frac{{\rm d}z_i}{z_i}
\end{equation}
on $(\CC^*)^r$. By Lemma~\ref{lem:isinvariant},
${\rm d}L_G(s,x)=\varphi^*(\omega_s)_x$. Consequently, for fixed $x$, the kernel of the linear map
$\mathcal K_G\rightarrow T_x^*\mathcal M_{0,n}$ given by $s\mapsto {\rm d}L_G(s,x)$ is precisely $\ker S_G(x)$. The rank of this map at a general point equals
$
\operatorname{rank}\mathcal S(G)-\kappa(G)=n-3
$
by~\eqref{eq:candidate}. It is also the generic rank of
${\rm d}\varphi$. Hence $\dim\mathcal M_G=n-3$,
and $\varphi$ is generically finite onto its image. Thus, for a general~$x$, 
\begin{equation}\label{eq: one-form vanish}
    s\in\ker S_G(x)
\quad\Longleftrightarrow\quad {\rm d}L_G(s,x)=0 \quad\Longleftrightarrow\quad
\omega_s\big|_{T_{\varphi(x)}\mathcal M_G}=0.
\end{equation}

Assume that the projection
$\mathcal V_G\to\mathcal K_G$ is not dominant. Then,
for general $s\in\mathcal K_G$, the form ${\rm d}L_G(s,{\bf \cdot})$ has no
zero on $\mathcal M_{0,n}$. Equivalently, $\omega_s$ has no zero on the smooth locus of $\varphi({\cal M}_{0,n})$. Furthermore, a general $\omega_s$ has no zero on the smooth locus of the boundary $B \coloneqq \mathcal M_G^{\rm reg} \setminus \varphi({\cal M}_{0,n})$. Indeed, apply the fiber-dimension theorem to the incidence correspondence 
$${\cal W}\coloneqq\{(z,s)\in B\times\mathcal{K}_G:\omega_s\big|_{T_{z}\mathcal M_G}=0\}.$$
We have $\dim({\cal W})\leq \dim B + r-(n-3) \le n-4 + r - (n-3) = r-1$, so the projection~${{\cal W} \to \mathcal{K}_G}$ is not dominant. By~\cite[Theorem~1.11]{EsterovGB}, there are a variety $Y\subseteq(\CC^*)^{r-1}$~and~a~homomorphism
$$
q:(\CC^*)^r\rightarrow(\CC^*)^{r-1} \quad \text{such that} \quad \mathcal M_G=q^{-1}(Y).
$$
The identity component of ${\rm ker} (q)$ is a
positive-dimensional subtorus. Choose a nonzero one-parameter subgroup of this torus, represented by
${(v_1,\ldots,v_r)\in\ZZ^r}$. Through every point
$(z_1,\ldots,z_r)\in\mathcal M_G$ there then passes the curve
$ \lambda\mapsto (z_1\lambda^{v_1},\ldots,z_r\lambda^{v_r})$,
whose tangent vector at $\lambda=1$ is $(z_1v_1,\ldots,z_rv_r)\in T_z\mathcal M_G$.
If $s=\sum_i t_i u_i\in\ker S_G(x)$, then equation
\eqref{eq: one-form vanish} implies
$$
0\,=\,\omega_s(z_1v_1,\ldots,z_rv_r)
  \,=\,\sum_{i=1}^r t_iv_i.
$$
Thus, for every general $x$, $ \ker S_G(x)\subseteq H_v$,
where
$
H_v=
\left\{
\sum_{i=1}^r t_i u_i\in\mathcal K_G:
\sum_{i=1}^r t_iv_i=0
\right\}
$ is a proper hyperplane in $\mathcal K_G$.
Therefore, we obtain a contradiction to the condition~\eqref{eq:fromxtou2}:
\[ \sum_{\substack{x\in\PP^{n-1}\setminus\mathcal H\\x\ {\rm general}}}
\ker S_G(x)
\,\,\subseteq \,\,H_v
\,\,\subsetneq \,\, \mathcal K_G.\]
Thus $\mathcal V_G\to\mathcal K_G$ is dominant. Since, by \eqref{eq:candidate},  
 $\dim \mathcal{V}_G = \dim \mathcal{K}_G$, this map is also finite-to-one.
 
It remains to show that the map to $\mathcal{M}_{0,n}$ is dominant. Equivalently, the fiber over  a generic point
 $x \in \PP^{n-1}$ is non-empty. This follows from~\eqref{eq:candidate}. The corresponding linear system should have a nonzero solution. Equivalently, $\mathrm{rank}(S_G(x)) < \min(|G|, 2n)$ for a generic $x \in \PP^{n-1}$. 
 The left-hand side is the rank of the scattering matroid $\mathcal{S}(G)$, and by Lemma~\ref{lem:scatmat} this rank is bounded by $2n-3$. Thus, it is enough to have $\mathrm{rank}(\mathcal{S}(G))< |G|$.
\end{proof}

We can use linear algebra over $\mathbb{Q}$ to verify
that a graph is matroidally copious.
This is clear for (\ref{eq:candidate}) and \eqref{eq:ScatMatrank}.
For the equivalent conditions~(\ref{eq:fromxtou2})–(\ref{eq:fromxtou1}), we
take a finite sample of points $x \in \PP^{n-1}\backslash \mathcal{H}$,
and check (\ref{eq:fromxtou2}) for the finite sum over those points or
(\ref{eq:fromxtou1}) for the finite intersection.
This also 
gives a probabilistic method for certifying that a graph is~copious.

\section{Commutative Algebra} \label{sec4}

In this section we study the vanishing ideal $I_G$ of the scattering correspondence
$\mathcal{V}_G$. We call this the {\em scattering ideal}.  It lives in
the polynomial ring $\CC[s,x]$ with variables $s_{ij}$
for the edges in $G$ and variables $x_1,\ldots,x_n$ for the
vertices of $G$. This is the homogeneous coordinate ring
of $\,\PP^{|G|-1} \times \PP^{n-1}$. The prime ideal $I_G$ is homogeneous in the
$\ZZ^2$-grading on $\CC[s,x]$ given~by
\begin{equation}
\label{eq:bigrading}
 {\rm deg}(s_{ij}) = (1,0) \quad {\rm and} \quad {\rm deg}(x_i) = (0,1). 
\end{equation}

Let $K_G$ denote the ideal generated by the $n$ linear forms in (\ref{eq:momcon})
which express momentum conservation. As a first approximation to $I_G$, we write
 $I^{(0)}_G$ for the ideal generated by $K_G$~and
 \begin{equation}
 \label{eq:naiveclearing}
 \prod_{j \in G_i} (x_i - x_j)  \cdot \frac{\partial L_G}{\partial x_i} 
\qquad {\rm for} \,\, i \in [n].
 \end{equation}

\begin{proposition} \label{prop:atmostone}
The ideal $I^{(0)}_G$ has at most one associated prime 
which contains no polynomial in $\CC[x]\backslash \{0\}$ and whose intersection 
with $\CC[s]$ equals $K_G$. That prime ideal equals~$I_G$.
\end{proposition}

\begin{proof}
The ideal $I_G$ is prime because it is the radical ideal of a vector bundle
over a Zariski open subset of $\PP^{n-1}$.
In the definition of $\mathcal{V}_G$, 
the point $s$ is generic in $\mathcal{K}_G$ and 
$x$ is generic in $ \PP^{n-1}\backslash \mathcal{H}$.
This implies that $I_G$ contains no
nonzero polynomial in $x$, and the only polynomials in $s$
that appear in $I_G$ are those in $K_G$. 
From polynomials in the ideal $I^{(0)}_G$, we therefore
remove any such extraneous factors in $x$ or $s$ alone.
After eliminating $\kappa(G)$ of the Mandelstam variables $s_{ij}$ from the linear ideal $K_G$, we can write our requirements as follows:
\begin{equation}
\label{eq:requirement}
 I_G \,\, = \bigl(S^{-1}I_G^{(0)}\bigr)\cap \CC[s,x] \,\, = \,\, \bigl(\CC(s)[x] \cdot I^{(0)}_G \,+ \, \CC(x)[s] \cdot I^{(0)}_G  \bigr) \,\cap \, \CC[s,x],
\end{equation} 
where $S = \{g(s)q(x): 0\neq g \in \CC[s] \, \text{ and } \,  0\neq q \in \CC[x]\}$. Note that the sum of the two extended ideals above is taken inside the smallest subring
of $\CC(s,x)$ containing both $\CC(s)[x]$ and $\CC(x)[s]$.
This ring is the image of the injective map $\CC(s)\otimes_{\CC}\CC(x)\hookrightarrow\CC(s,x)$ given by
$$
\frac{f(s)}{g(s)}\otimes\frac{p(x)}{q(x)}
\,\mapsto\,
\frac{f(s)p(x)}{g(s)q(x)}.
$$
Thus, $I_G$ is the intersection of all primary components of $I^{(0)}_G$ 
which contain no polynomial in $\CC[x]\backslash \{0\}$ and whose intersection 
with $\CC[s]$ is $K_G$. But, we already know
that $I_G$ is prime.
\end{proof}

\begin{definition} \label{def:algcop}
  A graph $G \subseteq \binom{[n]}{2}$ is \emph{algebraically copious} if
  the ideal $I_G$ defined in (\ref{eq:requirement}) 
   is a proper prime ideal in the ring $\CC[s,x]$, and its codimension is
$\operatorname{codim}(I_G)=\kappa(G)+n-3$.
\end{definition}

 Let $I^{(1)}_G$ be the ideal generated by $I^{(0)}_G$ and the polynomials in (\ref{eq:alwaysthere}).
This is a subideal of~$I_G$. Proposition \ref{prop:atmostone} holds verbatim for $I^{(1)}_G$.
In some cases, we are now done, i.e.~$I^{(1)}_G = I_G$.

\begin{example}  \label{ex:Kn} The complete graph $K_n$ is algebraically copious.
The ideal $I^{(0)}_{K_n}$ is generated by
$n$ polynomials in~\eqref{eq:momcon} of degree $(1,0)$
and $n$ polynomials in \eqref{eq:naiveclearing} of degree $(1,n\!-\!2)$.
The~ideal 
$$\,K_{K_n} + \langle \,x_i - x_j \,:\, 1 \leq i \!<\! j \leq n \rangle \,$$
is an embedded prime of $I^{(0)}_{K_n}$. Saturating by this embedded prime,
we arrive at $\,I_{K_n}^{(1)} = I_{K_n}$. Indeed, it follows from freeness of the braid arrangement \cite[Example 2.1]{KKMSW}
that $I_{K_n}$ is the complete intersection ideal generated by $K_{K_n}$
and polynomials (\ref{eq:alwaysthere}) for $k=2,3,\ldots,n-2$.
\end{example}

To illustrate Definition \ref{def:algcop}, we next show two graphs that are not
algebraically copious.

\begin{example}[$n=6$] \label{ex:witnesses}
Let $G_1$ be the graph with non-edges $\{35, 36, 45, 46 \}$,
and let $G_2$ be the graph with non-edges
$\{16,24, 25, 34, 35\}$.
Using {\tt Macaulay2}~\cite{Macaulay2Book} computations, we find~that
 \vspace{-0.18cm}
 $$ \begin{matrix}
 & (x_1-x_2)(x_3-x_5)(x_3-x_6)(x_4-x_5)(x_4-x_6)(s_{15}-s_{26}) \in I_{G_1}^{(1)} &{\rm but} &
 s_{15}-s_{26} \not\in K_{G_1},
  \smallskip \\
\! {\rm and} \!\! &
 (x_1-x_6)(x_2-x_4)(x_2-x_5)(x_3-x_4)(x_3-x_5)(s_{12}-s_{36}) \in I_{G_2}^{(1)} & {\rm but} &
 s_{12} - s_{36} \not\in K_{G_2}.
\end{matrix}
 \vspace{-0.18cm}
 $$
 This shows that  $K_{G_i}$ is strictly contained
in the intersection of $\CC(x)[s] \cdot I_{G_i}^{(1)}$ with $\CC[s]$.
The scattering equations have no solution in $\mathcal{M}_{0,6}$ for generic $s \in \mathcal{K}_{G_i}$,
thus
 $\mathcal{V}_{G_i} = \varnothing$ for $i=1,2$.
\end{example}

Following Theorem~\ref{lem:matroidally_geometrically},
we now come to the next step in our proof of Theorem \ref{conj:main}.
The next lemma completes the equivalence of Definitions \ref{def:geomcop},
\ref{def:matrcop} and \ref{def:algcop}.
Hence, from now on, we simply use the adjective {\em copious} for a
     graph that satisfies any of these three definitions.

\begin{lemma}\label{lem: geom. iff alg. copious}
A graph $G$ is geometrically copious if and only if it is algebraically copious.
\end{lemma}

\begin{proof}
We first prove the “only if” direction.
The projection $\mathcal{V}_G \to \mathbb{P}^{n-1}$ is dominant if and only if the
  elimination ideal $I_G \cap \mathbb{C}[x]$ is zero. Hence, the corresponding conditions in the definitions of geometrically and algebraically copious graphs are equivalent.
Now assume that the map $\mathcal{K}_G \leftarrow \mathcal{V}_G $ is finite-to-one modulo the
$\operatorname{PGL}(2)$ action. In particular, it is dominant onto $\mathcal{K}_G$, 
which is equivalent to $I_G \cap \mathbb{C}[s] = K_G$. Its generic fiber has dimension two, so $\operatorname{codim}(I_G)=\kappa(G)+n-3$.
Since \(I_G\) is a prime ideal, \(G\) is algebraically copious.

Conversely, if a graph  $G$ is algebraically copious, then $\mathcal{K}_G \leftarrow \mathcal{V}_G $  must also be finite-to-one. Indeed, the condition $I_G\cap\mathbb{C}[s]=K_G$ gives dominance over $\mathcal K_G$,
while the codimension condition gives two-dimensional generic fibers.
Hence the map is finite modulo $\operatorname{PGL}(2)$.
\end{proof}

We next introduce one more ideal for approximating $I_G$.
This is based on the connection with modules of derivations  on graphic arrangements \cite[Section~4]{KKMSW}.
The {\em derivation module} of $G$ is a $\CC[x]$-module. It consists of all vectors
$f(x) = (f_1(x),\ldots,f_n(x)) \in \CC[x]^n$ such that
\begin{equation}
\label{eq:isapolynomial} \qquad
\,\sum_{i=1}^n f_i(x) \, \frac{\partial L_G}{\partial x_i} \qquad \hbox{is a polynomial in $\,\CC[s,x]$}.
\end{equation}
Clearly, the polynomial in (\ref{eq:isapolynomial})
 is an element in our ideal $I_G$.
M\"uhlherr \cite[Theorem 1.3]{Muh} proved that the derivation module is generated by 
the finite set of {\em separator-based derivations}.
These are defined as follows.
    A subset $T\subset [n]$ is a \emph{separator} of the  graph $G$  if the induced subgraph $G \backslash T$ has more connected components than $G$. 
        Every connected component $C$ of $G\backslash T$ gives rise to a
          separator-based derivation $f(x)$. The $i$-th coordinate of this vector equals
\begin{equation}
\label{eq:f_i}  \qquad f_i(x) \,\,=\,\, \prod\limits_{t\in T} (x_i-x_t)  \,\,\,\,\text{ if } i\in C,  \quad 0\ \text{otherwise}.
\end{equation}
Then (\ref{eq:isapolynomial}) holds by \cite[Lemma 3.4]{Muh}, and
we obtain a polynomial of degree $(1,|T|-1)$ in $I_G$.

\begin{example}
Let $i$ be a vertex of $G$ that is not connected to every other vertex.
Then the set $G_i$ is a separator and $C = \{i\}$ is a connected component of $G\backslash G_i$.
Here (\ref{eq:isapolynomial}) is the generator
of $I_G^{(0)}$  shown in (\ref{eq:naiveclearing}). We also find this
   by clearing denominators in $\partial{L_G}/\partial x_i$.
\end{example}

We define $I^{(2)}_G$ to be the ideal in $\CC[s,x]$ that is generated by
$I^{(1)}_G$ and the polynomials (\ref{eq:isapolynomial}) for all
pairs $(C,T)$ as above. We summarize our discussion in the following corollary.

\begin{corollary} \label{cor:I2}
The ideal $I^{(2)}_G$ is contained in $I_G$. It also satisfies the conclusion of
Proposition~\ref{prop:atmostone}, namely:
if $\,I^{(2)}_G$ is prime, $\,I^{(2)}_G \cap \CC[x]  = \{0\}$ and
  $\,I^{(2)}_G \cap \CC[s] =  K_G$, then
$\,I^{(2)}_G = I_G$.
\end{corollary}

An important invariant of any $\ZZ^2$-homogeneous ideal  in
$\CC[s,x]$ is its {\em multidegree}. For the scattering ideal $I_G$ this is the class of $\overline{\mathcal{V}}_G$ in the cohomology ring of $ \PP^{|G|-1} \times \PP^{n-1}$. We have
\begin{equation}
\label{eq:multidegree} \qquad [\,\overline{\mathcal{V}}_G \,] \,\,  = \,\, \sum_{i=0}^c \gamma_i \,\sigma^i \tau^{c-i}
\,\, \,\,\in \,\,\ZZ[\sigma,\tau] /\langle \,\sigma^{|G|}, \,\tau^n \,\rangle\,\,
= \,\, H^*\bigl( \PP^{|G|-1} \times \PP^{n-1}, \ZZ\bigr).
\end{equation}
Here $c = {\rm codim}(\mathcal{V}_G)$, and
$\gamma_i = \# \bigl( \mathcal{V}_G \cap ( L_i \!\times\! M_{c-i}) \bigr)$,
where $L_i$ is a general subspace of dimension $i$ in $\PP^{|G|-1}$
and  $M_{c-i}$ is a subspace of dimension $c-i$ in $\PP^{n-1}$.
See Lemma~\ref{lem:multideg}.

\begin{proposition} \label{prop:acograph}
Let $G$ be a copious graph.
Then the sum  in (\ref{eq:multidegree}) for $[\overline{\mathcal{V}}_G]$
runs from $i=\kappa$ to $c = \kappa+n-3$,
where $\kappa=\kappa(G)$ is the codimension of the kinematic space  $\mathcal{K}_G$ in $\PP^{|G|-1}$. 
The last coefficient is $\gamma_{n+\kappa-3} = 1$, and the first nonzero coefficient gives the ML~degree:
\begin{equation}
\label{eq:muG} \mu(G) \,\, = \,\,\frac{ \gamma_\kappa}{n-2}. 
\end{equation}
\end{proposition}

\begin{proof}
The kinematic space $\mathcal{K}_G$ 
has dimension $|G|-\kappa-1$.
For copious graphs, the general
fiber of the map $ \mathcal{K}_G \leftarrow \mathcal{V}_G $ has dimension $2$.
This implies ${\rm dim}(\mathcal{V}_G) = |G|-\kappa+1$, and thus $c = n+\kappa-3$.
Moreover, that fiber is the union of $\mu(G)$ surfaces, each representing
one critical point of $L_G$ in $\mathcal{M}_{0,n}$. This union is a surface of degree
$\gamma_\kappa = \mu(G) \cdot (n-2)$. The rest of the statement is clear. Note that since $G$ is matroidally copious, we have $c < |G|$.
\end{proof}

 For instance, using the complete intersection in Example \ref{ex:Kn}, 
 the complete graph $K_n$ has
 $$ [\overline{\mathcal{V}}_{K_n}] \,\, = \,\, \sigma^n \prod_{k=2}^{n-2} (\sigma+ k \tau)  \,\, = \,\,
  (n-2)! \,\sigma^n \tau^{n-3} \,+\,\,\cdots \,\,+\,
 \frac{n(n-3)}{2}\, \sigma^{2n-4} \tau \,+ \,
 \sigma^{2n-3}  .
$$
From equation (\ref{eq:muG}) we now see the well-known formula
$ \,\mu(K_n) \,\,= \,\, (n-2)! /(n-2) \,\, = \,\, (n-3)!$.

\begin{example}[$n=5$] \label{ex:ideals5}
For each copious graph $G$,
we display the ideal~$I_G$ and its multidegree:
\begin{itemize}
\item[(a)] The complete graph $G = K_5$ has
$\, [\overline{\mathcal{V}}_{K_5}] = 6 \sigma^5 \tau^2 + 5 \sigma^6 \tau + \sigma^7 $, and its scattering
correspondence $\mathcal{V}_{K_5}$ is the complete intersection given by
(\ref{eq:alwaysthere}) for $k=2$ and $k=3$.
\item[(b)] If $G$ has one non-edge $45$, then
$\,[\overline{\mathcal{V}}_G] = 3 \sigma^5 \tau^2 + 4 \sigma^6 \tau + \sigma^7$.
The ideal  has three minimal generators, 
two of degree $(1,2)$ and one of degree $(2,2)$. We take
(\ref{eq:alwaysthere}) for $k=2$, and
$$ \begin{matrix} & & (x_1-x_2) (x_3-x_4) s_{24} \,\,+ \,\,(x_1-x_3) (x_2-x_4) s_{34}, \\ & {\rm and} &
(x_2-x_4) (x_3-x_5) s_{25} s_{34} \,-\, (x_2-x_5) (x_3-x_4) s_{24} s_{35}. \end{matrix}
$$
\item[(c)] If $G$ has two non-edges $23$ and $45$, then
$\,[\overline{\mathcal{V}}_G] = 3 \sigma^5 \tau^2 + 3 \sigma^6 \tau + \sigma^7$.
The ideal $I_G$ has three minimal generators modulo $K_G$, all  of degree $(1,2)$. We take
(\ref{eq:alwaysthere}) for $k=2$, and
$$  \begin{matrix} & & 
 (x_1-x_2) (x_3-x_5)   s_{25} \,+\, (x_1-x_3) (x_2-x_5)   s_{35}, \\ & {\rm and} & 
 (x_1-x_4) (x_3-x_5)   s_{34} \,+\, (x_1-x_5) (x_3-x_4)   s_{35}.
 \end{matrix}
$$
\end{itemize} 
We briefly examine Corollary \ref{cor:I2} for these
 graphs. We have
$I^{(2)}_G = I_G$ in (a) and~(c), but not in (b).
All generators of $I^{(2)}_G$ have degree $(1,\star)$,
but $I_G$ of (b) has a generator of degree~$(2,2)$.
\end{example}

We seek  $[\overline{\mathcal{V}}_G]$
for any  copious graph $G$. 
The coefficients  $\gamma_{i,G}$ can be 
found by solving polynomial equations as follows.
We pick $|G|-i-1$
random linear forms 
$\ell_1(s),\ldots,\ell_{|G|-i-1}(s)$
in the unknowns $s_{ij}$ for $ij \in G$.
These define a linear subspace
of dimension $i-\kappa$  in $\mathcal{K}_G$, where $\kappa = \kappa(G)$ from (\ref{eq:kappa}).
Next, we pick $i-\kappa+2$ random linear forms
$\ell'_1(x),\ell'_2(x),\ldots,\ell'_{i-\kappa+2}(x)$
in the unknowns $x_1,\ldots,x_n$.
These forms define a linear subspace of dimension
$\kappa + n-3-i$ in $\PP^{n-1}$.
The total number $|G|-\kappa+1$ of linear forms $\ell_i(s)$ and $\ell'_j(x)$
equals the dimension of $\mathcal{V}_G$ in $\PP^{|G|-1}\times \PP^{n-1}$.
As a result of this discussion, we have the following statement.
\begin{lemma} \label{lem:multideg}
The multidegree of the scattering correspondence is given by the formula
$$ \gamma_{i,G} \,\,=\,\, \# \,\bigl\{ \,(s,x) \in \mathcal{V}_G \,: \,
\ell_1(s) = \,\cdots \,= \ell_{|G|-i-1}(s) \,=\, \ell'_1(x)=\, \cdots \,=\ell'_{i-\kappa+2}(x)\, =\, 0\,
\bigr\} . $$
\end{lemma}

We now examine the case $n=6$.
The $2^{15}$ subsets of $\binom{[6]}{2}$ yield
 $156$ non-isomorphic graphs.
Most are not copious.
For instance, only $17$ of the $156$ graphs satisfy the
rank condition in \eqref{eq:candidate}.
This includes the two graphs
that were ruled out in Example~\ref{ex:witnesses}.
The following result summarizes our findings for $n=6$.
The cases $n=7,8$ will be treated in Section \ref{sec8}.

The proof of Proposition \ref{prop:14} is by direct computation,
mostly using  {\tt Macaulay2}~\cite{Macaulay2Book}.
The column for $[\overline{\mathcal{V}}_G]$ is found by the
command {\tt multidegree}.
From this we get $\mu(G)$ by (\ref{eq:muG}).
The rightmost column lists the $\ZZ^2$-degrees of the minimal generators of the
scattering ideal $I_G$. The index is the number of generators in that degree.
For instance, the symbol $(12)_3$ means that the ideal $I_G$ has three minimal generators
of degree $(1,2)$. The column $\delta(G)$ displays the degree of
the logarithmic discriminant, which will be discussed in Section \ref{sec9}.
The first graph in the list is the complete bipartite graph $K_{3,3}$.
This is special because $\kappa(K_{3,3})= 5$.

\begin{proposition} \label{prop:14}
There are $15$ copious graphs for $n=6$. Their ideals $I_G$ are as follows:
$$
\begin{matrix}
\hbox{non-edges of $G$} &\kappa(G) & \mu(G) &\delta(G) & \hbox{multidegree} \,[\overline{\mathcal{V}_G}] & \hbox{mingens of} \,\,I_G 
\smallskip \\
12, \! 13, \! 23, \! 45, \! 46, \! 56 & 5 & 1 & 0  & 4\sigma^5 \tau^3 + 7 \sigma^6 \tau^2 + 3 \sigma^7 \tau + \sigma^8  & (12)_5, (22) \\
12, 23, 34 &6& 1 & 0 & 4 \sigma^6 \tau^3 + 10 \sigma^7 \tau^2 + 6 \sigma^8 \tau + \sigma^9  & (12)_3, (22)_2, (32) \\
12, 13, 23 & 6& 1 & 0 & 4 \sigma^6 \tau^3 + 9 \sigma^7 \tau^2 +  6 \sigma^8 \tau + \sigma^9  & (12)_3, (22)_3 \\
12, 34, 45, 56 &6& 1 & 0 & 4 \sigma^6 \tau^3 + 10 \sigma^7 \tau^2 +  5 \sigma^8 \tau + \sigma^9  & 
\! (12)_3, \!(13), \!(22)_2, \! (32) \\
12, 13, 23, 45 &6& 1 & 0 & 4 \sigma^6 \tau^3 + 9 \sigma^7 \tau^2 +  5 \sigma^8 \tau + \sigma^9  & (12)_3, (13), (22)_3 \\
12, 23, 34, 45 &6& 1 & 0 & 4 \sigma^6 \tau^3 + 8 \sigma^7 \tau^2 +  5 \sigma^8 \tau + \sigma^9  & (12)_4, (22)_2 \\
12, 23, 34, 45, 56 &6& 1 & 0 & 4 \sigma^6 \tau^3 + 7 \sigma^7 \tau^2 + 4 \sigma^8 \tau + \sigma^9  & (12)_5 , (22)\\
12, 15, 23, 34, 45 &6& 1 &0  & 4 \sigma^6 \tau^3 + 6 \sigma^7 \tau^2 + 4 \sigma^8 \tau + \sigma^9  & (12)_6 \\

12, 23 &6& 2 &   4  &  8 \sigma^6 \tau^3 +  15 \sigma^7 \tau^2 + 7 \sigma^8 \tau + \sigma^9  & 
\! (12)_2, \!(13),\! (23), \!(42) \\
12, 34, 56 &6& 2 &6 & 8 \sigma^6 \tau^3 + 14  \sigma^7 \tau^2 + 6 \sigma^8 \tau + \sigma^9  & (12), (13)_4, (22)_2\\
12, 34, 45 &6& 2 & 4 & 8 \sigma^6 \tau^3 +  13 \sigma^7 \tau^2 + 6 \sigma^8 \tau + \sigma^9  & (12)_2, (13)_2, (32) \\
15,25,36,46 & 6 & 2 &4 & 8 \sigma^6 \tau^3 + 10 \sigma^7 \tau^2 + 5\sigma^8 \tau + \sigma^9 & (12)_3, (13) \\
12, 34 &6& 3 & 10 & 12 \sigma^6 \tau^3 + 17 \sigma^7 \tau^2 +  7 \sigma^8 \tau + \sigma^9  & (12), (13)_3, (34)_2  \\
   12 &6& 4 & 16 & 16 \sigma^6 \tau^3 +  21\sigma^7 \tau^2 + 8 \sigma^8 \tau + \sigma^9  &
    \!(12), \!(13)_2, \!(35),\!(54) \\
\varnothing &6& 6 &30& 24 \sigma^6 \tau^3 + 26 \sigma^7 \tau^2 + 9 \sigma^8 \tau + \sigma^9 & (12), (13),(14) \\
\end{matrix}
$$
\end{proposition}

For readers of \cite{KKMSW} we note that
$12$ of the $15$ copious graphs are not gentle.
 Their ideals $I_G$ have minimal generators of
some degree $(2,\star)$ or $(3, \star)$.
Here, the primality criterion in Corollary \ref{cor:I2} does not apply,
and it requires a non-trivial saturation step to compute $I_G$ from its subideal $I_G^{(2)}$.
Adding to Example~\ref{ex:witnesses}, here is another
 interesting non-copious graph.

\begin{example}
Let $G$ be the edge graph of a triangular prism, with non-edges
 $12,  23, 34, 45, 56 , 16$.
The ideal $I^{(2)}_G$  is generated modulo $K_G$ by seven
polynomials of degree $(1,2)$. It is radical and has two associated primes.
One of these primes is
the irrelevant ideal  $\langle s_{13}, s_{14}, s_{15}, s_{24},s_{25},s_{26}, s_{35},s_{36},s_{46} \rangle$
of the projective space $\mathcal{K}_G \simeq \PP^8$.
The other associated prime has the multidegree
$4 \sigma^6 \tau^3 + 5 \sigma^7 \tau^2 + 3 \sigma^8 \tau$.
The absence of the term $\sigma^9$ shows that this prime ideal contains a
polynomial in $\CC[x] \backslash \{0\}$. Indeed, we found the cubic generator
\begin{equation}
\label{eq:cubicinx} \begin{matrix}
\quad x_1 x_2 x_4 - x_1 x_2 x_5 - x_1 x_3 x_4 + x_1 x_3 x_6 + x_1 x_4 x_5 - x_1 x_4 x_6 \\
 + \,x_2 x_3 x_5 - x_2 x_3 x_6 - x_2 x_4 x_5 + x_2 x_5 x_6 + x_3 x_4 x_6 - x_3 x_5 x_6. \,
 \end{matrix}
 \end{equation}
 Vanishing of this cubic is a necessary condition for a point in $\mathcal{M}_{0,6}$ to
 be a critical point of some scattering potential on $G$.
 We conclude that $\mathcal{V}_G = \varnothing$, meaning that $G$ is not copious.
  \end{example}

One may wonder whether the class of copious graphs is
closed under adding edges.
The graph $K_{3,3}$ in the first row of 
Proposition  \ref{prop:14} shows that this is not the case.
Let $G'$ be obtained  from $K_{3,3}$ by adding one edge $ij$. 
The scattering matroid $\mathcal{S}(K_{3,3})$ has rank $8$, but
 $\mathcal{S}(G')$ has rank $9$. The new edge $ij\in G'$ is a co-loop,
 so condition (\ref{eq:candidate}) is violated.
 Namely, the kernel of each $S_G(x)$ is contained in a proper subspace of $\mathcal{K}_G$ given by $s_{ij}=0$.
 However, if we disallow bipartite graphs then the closure property holds.
  We state this result below. 
\begin{proposition}\label{prop: add edges}
 Let $G$ be a copious graph which is connected and not bipartite, 
 and let $G'$ be a graph obtained from $G$ by adding one edge.
 Then the larger graph $G'$ is also copious.
 \end{proposition}
 \begin{proof}
  Let $e$ be the new edge. Since $G$ and $G'$ are connected and
non-bipartite, $\kappa(G)=\kappa(G')=n$. By copiousness of $G$,
$\operatorname{rank}\mathcal S(G)=\operatorname{rank}\mathcal S(G\setminus f)=2n-3$ for each edge $f\in G$. By Lemma~\ref{lem:scatmat}, this is the maximal possible rank of a scattering matroid. Hence the same rank equalities hold for $G'$
and all its edge deletions $G'\setminus f$. It remains to verify condition~\eqref{eq:fromxtou2}.
We have $\mathcal K_G=\mathcal K_{G'}\cap\{s_e=0\}$. For general $x$, setting $s_e=0$ identifies
$\ker S_G(x)$ with a subspace of $\ker S_{G'}(x)$ of codimension one, so $\ker S_{G'}(x)$ contains some $v$ with
$v_e\neq0$. Hence, using condition~\eqref{eq:fromxtou2} for $G$, we have
$\,\sum_x\ker S_{G'}(x) \,\supseteq \,\mathcal K_G+\CC v 
\,=\, \mathcal K_{G'}$.
The reverse inclusion is automatic. This proves that $G'$ is matroidally copious, and  therefore
copious.
\end{proof}

\section{Partial Compactifications}\label{sec5}

We now turn to  
Huh's Theorem  \cite{Huh}, which equates  the 
signed Euler characteristic of a smooth very affine variety $X$ with its ML degree ${\rm MLdeg}(X)$.
The very affine variety in our setting is a  partial compactification
$\mathcal{M}_G$ of the moduli space $\mathcal{M}_{0,n}$.
This already partially appeared in Lemma  \ref{lem:isinvariant} and in the proof of Theorem~\ref{lem:matroidally_geometrically}. We will now recall the definition.

We follow the proof of Lemma  \ref{lem:isinvariant}.
The {\em kinematic lattice} $  \mathcal{K}_G \,\cap \, \ZZ^{|G|}$ is a 
free abelian group of rank $r = |G|-\kappa(G)$. 
Every vector $u \in  \mathcal{K}_G \,\cap \, \ZZ^{|G|}$
gives rise to a character $\tilde u : (\CC^*)^{|G|} \rightarrow \CC^*$.
We write $\tilde u$ as a Laurent monomial  in the
 $x_i - x_j$ for $ij \in G$. The generalized cross-ratio $\tilde u$ is constant
along ${\rm PGL}(2)$ orbits  on $\mathbb{C}^{n}$, so it is a regular function on
 the moduli space~$\mathcal{M}_{0,n}$.

Let $u_1,\ldots,u_r$ be a basis of 
 $  \mathcal{K}_G \cap  \ZZ^{|G|}$. The generalized cross-ratios define a map
 \begin{equation}
 \label{eq:umap1}
   \PP^{n-1} \,\dashrightarrow \, (\CC^*)^r, \,\,
 x \mapsto \bigl(\tilde u_1(x),\tilde u_2(x),\ldots,\tilde u_r(x) \bigr). 
 \end{equation}
 This rational map descends to a well-defined regular map on the moduli space:
  \begin{equation}
 \label{eq:umap2}
  \mathcal{M}_{0,n} \,\rightarrow \,(\CC^*)^r, \,\,
 x \mapsto \bigl(\tilde u_1(x),\tilde u_2(x),\ldots,\tilde u_r(x) \bigr). 
\end{equation}  
We define $\mathcal{M}_G$ as the Zariski closure in $(\CC^*)^r$
of the image of (\ref{eq:umap1}) or (\ref{eq:umap2}). Then
$\mathcal{M}_G$ is a very affine variety. This means that
$\mathcal{M}_G$ is the variety defined by an ideal $\mathcal{I}_G$ in the
Laurent polynomial ring $\CC[z_1^{\pm 1} , z_2^{\pm 1}, \ldots, z_r^{\pm 1}]$,
whose elements are the regular functions on $(\CC^*)^r$.
We note that $\mathcal{M}_G$ and $\mathcal{I}_G$ do not
really depend on the choice of basis for the kinematic lattice
$  \mathcal{K}_G \,\cap \, \ZZ^{|G|}$. Passing to a different basis
induces a  change of coordinates on the torus $(\CC^*)^r$.

\begin{example}[Octahedron]  \label{ex:octah3}
Fix 
$G$ in Example \ref{ex:octah} and the basis
(\ref{eq:octabasis1}) for its kinematic lattice.
The map $\PP^{5} \dashrightarrow (\CC^*)^6$  in (\ref{eq:umap1}) has the six coordinates
$z_i = \tilde u_i(x)$ from
(\ref{eq:octabasis2}).
Elimination reveals that the very affine variety $\mathcal{M}_G$ is the smooth complete intersection~given~by
$$ \mathcal{I}_G \,\,=\,\,  \bigl\langle
z_1 z_4 z_5-z_4-z_5+1, \,z_2 z_4 z_6-z_2 z_4-z_2 z_6+1, \,z_3 z_5 z_6-z_3 z_5+z_3-z_6
\bigr\rangle .
$$
Assuming momentum conservation (\ref{eq:momcon}),
the scattering potential (\ref{eq:octaLG}) can be written as
$$ L_G(t,z) \, = \, 
t_1 \,{\rm log}(z_1) \,+\,
t_2 \,{\rm log}(z_2) \,+\,
t_3 \,{\rm log}(z_3) \,+\,
t_4 \,{\rm log}(z_4) \,+\,
t_5 \,{\rm log}(z_5) \,+\,
t_6 \,{\rm log}(z_6) ,
$$
where the coefficients $t_i$ are  linear combinations of the $s_{ij}$.
We restrict this function to  the subvariety $\mathcal{M}_G$ of $(\CC^*)^6$.
Using Lagrange multipliers for $\mathcal{I}_G$, we find that
this restriction has $\mu(G) = 2$ critical points.   Thus the
threefold $\mathcal{M}_G$ has Euler characteristic
$\chi(\mathcal{M}_G) = -2$. 
\end{example}

We now come to our fourth and last definition 
of what it means for a graph to be copious.

\begin{definition} \label{def:topcop}
  A graph $G \subseteq {[n] \choose 2}$ is \emph{topologically copious} 
if the map (\ref{eq:umap2}) is an embedding, that is, an isomorphism onto a locally closed subvariety of $(\mathbb C^*)^r$ for $r = |G| - \kappa(G)$.
This property is independent of the choice of basis
for the kinematic lattice $ \mathcal{K}_G \,\cap \, \ZZ^{|G|}$.
\end{definition}  

\begin{theorem} \label{thm:pacom}
For a topologically copious graph $G$,
 the very affine variety $\mathcal{M}_G$
is a partial compactification of
$\mathcal{M}_{0,n}$. Moreover,
$ \mu(G)=\operatorname{MLdeg}(\mathcal M_G)$, and if
$\mathcal M_G$ is smooth, these are also equal to the signed
topological Euler characteristic $(-1)^{n-3}\chi(\mathcal M_G)$.
\end{theorem}

\begin{proof}[Proof and Discussion]
The map (\ref{eq:umap2}) is an embedding of
$\mathcal{M}_{0,n}$ into the torus $(\CC^*)^r$.
The very affine variety $\mathcal{M}_G$ is the closure
of $\mathcal{M}_{0,n}$ in $(\CC^*)^r$.
This closure is a partial compactification of $\mathcal{M}_{0,n}$ and hence has dimension $n-3$.
For a similar construction see~\cite{Fry2022}. Suppose that $\mathcal{M}_G$ is smooth. By Huh's Theorem~\cite{Huh}, its signed Euler characteristic $(-1)^{n-3}\chi(\mathcal M_G)$ equals the ML degree ${\rm MLdeg} ({\cal M}_G)$, which is the number of critical points of a generic log-likelihood function. 
That function is our scattering potential $L_G$, as seen from 
(\ref{eq:LGinvariant}). 
Since the boundary
$\mathcal M_G\setminus\mathcal M_{0,n}$ has smaller dimension, the critical-point incidence over it does not dominate
$\mathcal K_G$. Hence, for generic $s\in\mathcal K_G$, all
critical points lie in $\mathcal M_{0,n}$.
This proves
the last assertion. 
\end{proof}

Graphs $G$ for which $\mathcal{M}_G$ is not smooth appear in Section \ref{sec6}.
We now turn to a  special class of graphs $G$
for which $\mathcal{M}_G$ is the complement of
a hyperplane arrangement. Here, the very affine variety
$\mathcal{M}_G \subset (\CC^*)^r$ is cut out by linear equations,
so $\mathcal{M}_G$ is smooth.
A vertex of $G$ is a {\em universal vertex} if
it is connected to all other vertices in $G$.
This is the assumption made tacitly in many sources, including \cite{EPS, Lam, RS}. We now state our result
for this case:

\begin{proposition} \label{prop:gauge}
Every copious graph that has
a universal vertex
is topologically copious.
\end{proposition}

\begin{proof}
We write $[ij] = x_i - x_j$ for the $2 \times 2$ minors of our matrix
which coordinatizes $\PP^{n-1}$:
$$ X \,\, = \,\, \begin{pmatrix}
1 & 1 & 1 & \cdots & 1 \\
x_1 & x_2 & x_3 & \cdots & x_n 
\end{pmatrix}.
$$
A projective isomorphism yields
   the coordinate system on $\mathcal{M}_{0,n}$ referred to as {\em gauge fixing}:
\begin{equation}
\label{eq:gaugefixing}  Z \,\, = \,\, \begin{pmatrix}
1 & 1 & 1 &   1 & \cdots & 1 & 0 \\
0 & 1 & z_3 & z_4 & \cdots & z_{n-1} & 1 
\end{pmatrix}.
\end{equation}
The $2 \times 2$ minors of $Z$ define an arrangement $\mathcal{A}$ of
$\binom{n-1}{2} -1 $ affine hyperplanes in $ \CC^{n-3}$.
This is the coordinate choice used, e.g., in~\cite[equation (1)]{EPS}, 
\cite[Section 3]{KKMSW}, and
\cite[equation (2)]{ST}.

We can write the affine-linear forms for these hyperplanes as generalized cross-ratios:
\begin{equation}
\label{eq:linearcrossratios}
 \frac{[1i][2n]}{[12][in]} = z_i \,,\quad
 \frac{[1n][2i]}{[12][in]} = z_i - 1\, , \quad
 \frac{ [1n][2n][ij]}{[12][in][jn]} \, = z_j - z_i 
 \quad {\rm for} \,\, i,j \in \{3,4,\ldots,n-1\}. 
 \end{equation}
The corresponding elements of $\ZZ^{\binom{n}{2}}$ form a $\ZZ$-basis of 
the kinematic lattice $\mathcal{K}_{K_n} \cap \ZZ^{\binom{n}{2}}$:
\begin{equation}
\label{eq:latticebasis}
e_{1i} + e_{2n} - e_{12} - e_{in}\,,\,\,\,
e_{1n} + e_{2i} - e_{12}- e_{in}\,,\,\,\,
e_{1n}+e_{2n} + e_{ij} - e_{12} - e_{in} - e_{jn} .
\end{equation}

We now assume that $n$ is a universal vertex of $G$
and that $12$ is an edge of $G$.
Our choice of $Z$ ensures that $[12] = 1$ and $[in] = 1$ for all $i$.
Let $\mathcal{A}_G$ denote the subarrangement of $\mathcal{A}$
obtained by deleting all hyperplanes whose label $1i$,
$2i$ or $ij$ is not an edge of $G$.
The remaining linear forms in (\ref{eq:latticebasis}) form
a $\ZZ$-basis of the kinematic space $\mathcal{K}_G$.
The cross-ratios (\ref{eq:linearcrossratios}) for $ij \in G$ 
are regular functions on $\mathcal{M}_G$, written in the coordinates $z_i$ of $\CC^{n-3}$.

We identified
the very affine variety $\mathcal{M}_G$ with the arrangement complement
$\CC^{n-3} \backslash \mathcal{A}_G$.
Since $G$ is  copious, the map
$\mathcal{V}_G \rightarrow \mathcal{K}_G$ is
$\mu(G)$-to-1. In particular, $\mu(G) \geq 1$. By Varchenko's Theorem,
the ML degree $ \mu(G)={\rm MLdeg}({\cal M}_G)$ is the
number of bounded regions of $\mathbb{R}^{n-3}\setminus\mathcal{A}_G$.

We fix one bounded region of $\mathbb{R}^{n-3}\setminus\mathcal{A}_G$, and 
we consider the affine-linear forms
$z_i$, $z_i-1$  or $z_j-z_i$ which
define the facets of that region.
These span the space of all affine-linear forms on $\RR^{n-3}$.
In fact, since graphic matroids are unimodular, 
they contain a $\ZZ$-basis for the lattice of affine-linear
forms on $\ZZ^{n-3}$.
Hence, every coordinate $z_3,z_4,\ldots,z_{n-1}$
can be written as a polynomial in the generalized cross-ratios
(\ref{eq:linearcrossratios}) that come
 from our $\ZZ$-basis of $\mathcal{K}_G \cap \ZZ^{|G|}$.
This shows that the map 
$\mathcal{M}_{0,n} \rightarrow \mathcal{M}_G$ is an embedding, 
so
 $G$ is topologically~copious.
\end{proof}

From  the proof above, we can easily derive
 linear generators for the ideal $\mathcal{I}_G$
that defines $\mathcal{M}_G$ as a closed subvariety
in $(\CC^*)^r$. Namely, we consider the 
$r$ expressions (\ref{eq:linearcrossratios})
which correspond to edges of $G$.
These $r$ expressions are linear in the $z_i$.
The relations among them generate $\mathcal{I}_G$,
and these are linear as well.
If $G$ does not have a universal vertex then 
$\mathcal{I}_G$ requires  non-linear Laurent polynomials.
We saw this for the octahedron in Example  \ref{ex:octah3}.

\begin{remark}\label{rem:deconing}
    The affine arrangement $\mathcal{A}$ is the \emph{deconing}, in the sense of
    \cite{OrlikTerao2001, Stanley}, of the braid arrangement $\{x_i = x_j\}_{ij \in K_{n-1}}$
    with respect to one hyperplane. Note that this has $n-1$ variables, not $n$.
    For example, for $n = 5$, if we change coordinates to
    $z_0 =  x_2 - x_1$,
    $z_1 =  x_3 - x_1$
    and 
    $z_2 = x_4-x_1$, then
    intersecting with the plane     $z_0=1$  in $\RR^3$
    gives the line arrangement
        $\mathcal{A}$ in $\RR^2$.
    In general, $\mathcal{A}_G$ is the deconing of the graphical arrangement~$\mathcal{H}_{G \backslash n}$.
    \end{remark}
We now conclude the proof of our main theorem.
We also record a converse under the
additional assumption that the very affine variety $\mathcal M_G$ has positive ML degree $\mu(G)$.
\begin{proof}[Proof of Theorem~\ref{conj:main}]
The equivalence of geometric, matroidal, and algebraic copiousness
was established in Theorem~\ref{lem:matroidally_geometrically}
and Lemma~\ref{lem: geom. iff alg. copious}. If $G$ has a universal vertex,
Proposition~\ref{prop:gauge} shows that copiousness implies
topological copiousness and hence completes the proof.
\end{proof}
\begin{proposition}\label{prop: top copious + ML deg positive}
If $G$ is topologically copious and
$\mu(G)\ge 1$, then $G$ is copious.
\end{proposition}

\begin{proof}
Suppose that $G$ is topologically copious,
and consider the very affine variety $\mathcal{M}_G$
in~$(\CC^*)^r$. The results in \cite{Huh} ensure that the associated likelihood correspondence
$$
 \PP^{r-1} \,\leftarrow\, \mathcal{V}_G\, \rightarrow \,\mathcal{M}_{G}
$$
 has the desired properties, namely the right map is dominant and the left map is 
 $\mu(G)$-to-$1$. Here, $\mu(G) =\operatorname{MLdeg}(\mathcal M_G)\ge 1$. The inequality is our assumption.
  The same holds for the maps in (\ref{eq:twoprojections}). Indeed, as in proof of Theorem~\ref{thm:pacom}, generic critical points lie in $\mathcal M_{0,n}$. Moreover,  $\mathcal{K}_G = {\rm span}(u_1,\ldots,u_r) = \PP^{r-1}$.
This means that $G$ is geometrically~copious.
\end{proof}

\section{Complete Bipartite Graphs}\label{sec6}

In this section we present a case study for the
complete bipartite graph $G = K_{s,t}$.
We assume that $3\leq s\leq t $, so the number of vertices of $G$
satisfies $n = s+t \geq 6$. The $st$ edges of $G$ are 
the pairs $ij$, where $i \in \{1,2,\ldots,s\}$ and $j \in \{s+1,\ldots,n\}$.
The smallest instance $K_{3,3}$ has 
 $\kappa(K_{3,3}) = 5$ and its ML degree is $\mu(K_{3,3}) = 1$.
This is the first row in Proposition \ref{prop:14}.

\begin{lemma}
The kinematic space $\mathcal K_{K_{s,t}}$ has dimension
$(s-1)(t-1)$. A $\mathbb Z$-basis of 
$\mathcal K_{K_{s,t}}\cap\mathbb Z^{st}$
is given by the vectors $ u_{ij} = e_{1,s+1}+e_{ij}-e_{1j}-e_{i,s+1}$,
where $2\leq i\leq s$ and $s+2\leq j\leq n$.
\end{lemma}

\begin{proof}
Corollary \ref{cor: bipartite connected components} gives $\kappa(K_{s,t}) = s+t-1$ and hence $\dim\mathcal K_{K_{s,t}}=st-\kappa(K_{s,t})=(s-1)(t-1)$. Since the graph $K_{s,t}$ is bipartite, its even-cycle matroid is the graphic
matroid of $K_{s,t}$, and thus the kinematic space~$\mathcal{K}_{K_{s,t}}$
is the cycle space of $K_{s,t}$. This cycle space
is spanned by the $4$-cycles, and the listed $4$-cycles $u_{ij}$
form an integer basis for the kinematic lattice.
\end{proof}

The $(s-1)(t-1)$ cross-ratios defined by our cycle basis for $K_{s,t}$ are
\begin{equation}
\label{eq:crossratios44}
 \quad \tilde u_{ij}(x) \quad = \quad
\frac{[1,s\!+\!1] \,[ij]}{[1j] \,[i,s\!+\!1]} \,\, = \,\,
 \frac{(x_1 - x_{s+1})(x_i - x_j)}{(x_1 - x_j)(x_i- x_{s+1})}
\quad \qquad \begin{matrix} {\rm for} \,\,\,\,  i = 2,3,\ldots,s
\\  \; \,{\rm and} \,\, \ j = s\!+\!2, \ldots,n. \quad
\end{matrix}  
\end{equation}
The brackets $[ij]$ are  Pl\"ucker coordinates
for the Grassmannian ${\rm Gr}(2,n)$ from
which the moduli space $\mathcal{M}_{0,n}$ arises
by taking the quotient modulo $(\CC^*)^n$. We consider the matrix
\begin{equation}
\label{eq:Zminus1}  \mathcal{Z} \,\, \coloneqq \,\,
\begin{pmatrix}
z_{2,s+2}-1 & z_{2,s+3}-1 & \cdots & z_{2n} -1\\
z_{3,s+2}-1 & z_{3,s+3}-1 & \cdots & z_{3n} -1\\
 \vdots & \vdots & \ddots  & \vdots \\
z_{s,s+2}-1 & z_{s,s+3}-1 & \cdots & z_{sn} -1\\
\end{pmatrix}.
\end{equation}
This matrix has $s-1$ rows and $t-1$ columns.
They are indexed by $i$
and $j$ as in (\ref{eq:crossratios44}).
The entries $z_{ij}$ are variables, and we 
consider the map (\ref{eq:umap2}) that is
given by $z_{ij} = \tilde u_{ij}(x)$.

\begin{theorem}\label{thm:bipartite}
The ideal $\mathcal{I}_{K_{s,t}}$ is generated by the $2\times2$ minors of the
matrix ${\cal Z}$. The partial compactification $\mathcal M_{K_{s,t}}$
is the subvariety of $(\CC^*)^{(s-1)(t-1)}$ consisting of all
points~$(z_{ij})$ with $\mathrm{rank}({\cal Z})\leq 1$.
The map~\eqref{eq:umap2} is an embedding,
so the graph $K_{s,t}$ is topologically copious. The very affine variety
$\mathcal M_{K_{s,t}}$ has precisely one singular point, given by
$z_{ij}=1$ for all $i,j$.
\end{theorem}
\begin{proof}
We fix the gauge $x_1=\infty$, $x_{s+1}=0$, $x_{s+2}=1$.
The cross-ratios in~\eqref{eq:crossratios44} then become
$$ \begin{matrix} z_{ij}\,=\,1-\frac{x_j}{x_i}, \quad \text{so } \; z_{ij}-1\,=\,-\frac{x_j}{x_i}. \end{matrix} $$
Thus the matrix ${\cal Z}$ has rank at most one and $\mathcal M_G$ is contained in the variety
$\{\operatorname{rank}({\cal Z})\leq1\}$.

Conversely, let $U$ be the open subset of the variety
$\{\operatorname{rank}({\cal Z})\leq1\}$ on which all entries of ${\cal Z}$ are nonzero and the rows and columns of ${\cal Z}$ are pairwise distinct. On $U$, we have
$$ \begin{matrix}
x_i\,=\,-\frac{1}{z_{i,s+2}-1} \quad {\rm and}
\quad
x_j\,=\,\frac{z_{2j}-1}{z_{2,s+2}-1}.
\end{matrix}
$$
The $2\times2$ minors of ${\cal Z}$ show that these formulas recover all
the coordinates $z_{ij}$. They define a regular map $U\to\mathcal M_{0,n}$ which
is inverse to~\eqref{eq:umap2}, so the map~\eqref{eq:umap2} is an embedding.
The variety $\{\operatorname{rank}({\cal Z})\leq1\}$ is irreducible of
dimension
$
(s-1)+(t-1)-1=n-3,
$
and its~prime ideal is generated by the $2\times2$ minors of ${\cal Z}$.
The image is a dense open subset, so its closure is this
determinantal variety and its ideal is $\mathcal{I}_{K_{s,t}}$. Finally, the zero
matrix is the unique~singular point of the determinantal variety.
In our coordinates, this is the point $z_{ij}=1$ for~all~$i,j$.
\end{proof}

\begin{remark}
A general point of $\mathcal M_{K_{s,t}}$ admits a nice
matrix completion interpretation. We augment the rank one matrix
${\cal Z}$  by a first column of ones
and a first row $(1,0,0,\ldots,0)$, and we then
add the first column to each of the other columns.
The resulting $s \times t$ matrix $X$ has rank two, by construction. By \cite[Proposition~7.1]{SU}, we can find matrices
$P$ and $Q$ such that $\begin{footnotesize} \begin{pmatrix} P & X \\ \! -X^\top & Q \end{pmatrix}
\end{footnotesize}$ is skew-symmetric of~rank~2. 
The blocks $P$ and $Q$ recover the missing Pl\"ucker coordinates within the two parts of~$K_{s,t}$. This completion is unique up to the torus action $(P,Q)\mapsto(t P,t^{-1}Q)$.
After quotienting by this $\CC^*$ action, the matrix completion determines a
general point in $\mathcal M_{0,n}$. This elucidates 
the birationality of the cross-ratio~map~\eqref{eq:umap2}.
\end{remark}

Theorem~\ref{thm:bipartite} allows us to determine 
$\mu(K_{s,t})$ from results
in the algebraic statistics  literature.
Let $\mathcal{X}_{s,t}$ denote the 
 projective variety of 
 $s \times t$ matrices $X$ of rank $\leq 2$.
The ideal of $\mathcal{X}_{s,t}$ is generated
by the $3 \times 3$ minors. The positive points on
$\mathcal{X}_{s,t}$ form a  statistical model
inside the probability simplex $\Delta_{st-1} \subset \PP^{st-1}$.
This represents  conditional
independence  for two random variables,
with $s$ states and $t$ states, with binary 
latent random variable.
The algebraic study of likelihood inference for this model was initiated
by Hauenstein, Rodriguez and Sturmfels~\cite{HRS}.
They conjectured a formula for the ML degree, and that
conjecture was proved by Rodriguez and Wang in \cite{RW}.
Their topological approach also solves our problem.

\begin{theorem} \label{thm:jose}
The ML degree $\mu(K_{s,t})$ of the complete bipartite graph 
$K_{s,t}$ is equal to the ML degree of
the conditional independence model $\mathcal{X}_{s-1,t-1}$.
We have $\mu(K_{3,t}) = 1$, $\mu(K_{4,t}) = 2^t-6$,
and $\mu(K_{5,t} )= 2 \cdot 3^t - 5 \cdot 2^{t+1} + 25$.
Further values can be found in \cite[Table 1, page 501]{RW}.
\end{theorem}

\begin{proof}
Set $m=s-1$ and $q=t-1$, so $m\leq q$. By
Theorem~\ref{thm:bipartite}, the variety $\mathcal M_{K_{s,t}}$ equals
$${\cal Y}_{m,q} \,\,\coloneqq\,\, \left\{
Z\in(\CC^*)^{mq}:
\operatorname{rank}(Z-\mathbf 1)\leq1
\right\}.$$
Let $A=(a_{ij})$ be generic data and consider the likelihood function
${L_A(Z)\,\,=\,\,\sum_{i,j}a_{ij}\log(z_{ij})}$. For a  matrix $Q=(q_{ij})$, set
$q_{++}=\sum_{i,j}q_{ij}$ and $a_{++}=\sum_{i,j}a_{ij}$ and consider the function
$$L_A^{\rm proj}(Q) = \sum_{i,j}a_{ij}\log(q_{ij}) -a_{++}\log(q_{++}).$$
Its differential is  given by the matrix 
$
\left(
\frac{a_{ij}}{q_{ij}}-\frac{a_{++}}{q_{++}}
\right)_{ij}$.
At a smooth point of ${\cal Y}_{m,q}$, write ${Z=\mathbf1+xy^\top}$ and define $Q=(q_{ij})$ by $q_{ij}\coloneqq \frac{a_{ij}}{z_{ij}}$.
The tangent space ${\cal Y}_{m,q}$ at smooth $Z$ equals
$$
T_Z ({\cal Y}_{m,q})
=
\left\{v y^\top+x w^\top:
v\in\CC^m,\ w\in\CC^q\right\}.
$$
We have $ {\rm d}L_A(Z)(vy^\top+xw^\top)  =v^\top Qy+x^\top Qw$.
Thus $Z$ is a critical point of $L_A$ if and only if
$Qy=0$ and $x^\top Q=0$.
This implies $\operatorname{rank}(Q)\leq m-1$.
Moreover, $ a_{++} = \sum_{i,j}q_{ij}(1+x_i y_j) =
q_{++}+x^\top Qy= q_{++}$.  Consequently, the differential of $L_A^{\rm proj}$ at $Q$ is
represented by the matrix
$$
\left(\frac{a_{ij}}{q_{ij}}-1\right)_{ij}
=Z-\mathbf 1
=xy^\top.
$$

Since $Qy=0$ and $x^\top Q=0$, the matrix $xy^\top$ belongs to the normal space at $Q$ of the~variety
${\cal Q}_{m,q} \coloneqq \{\mathrm{rank}(Q) \le m-1\}$. 
Hence $Q$ is a likelihood critical~point~on~${\cal Q}_{m,q}$.

Conversely, the same calculation shows that every generic likelihood critical point $Q$ on ${\cal Q}_{m,q}$, normalized by $q_{++}=a_{++}$, gives a critical point $Z=(\frac{a_{ij}}{q_{ij}})=\mathbf1+xy^\top$ of ${\cal Y}_{m,q}$. Hence entrywise division gives a bijection between
the two generic critical loci. Therefore, we have
$$
\operatorname{MLdeg}({\cal Y}_{m,q})
\,\,=\,\,
\operatorname{MLdeg}({\cal Q}_{m,q})
\,\,=\,\,\operatorname{MLdeg}({\cal X}_{m,q}),
$$
where the last equality follows from Maximum Likelihood duality for determinantal varieties, see
\cite[Theorem~1]{DraismaRodriguez}. Thus, by Theorem~\ref{thm:pacom},
$\mu(K_{s,t}) = \operatorname{MLdeg}(\mathcal M_{K_{s,t}})
= \operatorname{MLdeg}(\mathcal X_{s-1,t-1})$.
In particular, since this number is positive,  Proposition~\ref{prop: top copious + ML deg positive} implies that the graph
$K_{s,t}$ is copious. The desired
formulas then follow from~\cite[Theorem~3.12 and Corollary~4.2]{RW}.
\end{proof}

\begin{example}
Consider the graph $K_{3,t}$. By Theorem~\ref{thm:jose}, its ML
degree is one. This means that the Maximum Likelihood Estimator (MLE) can be written as a rational function in the data, and
every coordinate of that rational function has a Horn factorization \cite{EPS}
as an alternating product of positive linear forms.
In what follows, we present the formula.

The likelihood problem is to find the critical points of
$\sum_{i=1}^2\sum_{j=1}^{t-1}a_{ij}\log(z_{ij})$
on the variety of rank one matrices
$ Z-{\bf 1}=(z_{ij}-1)$. For generic data, the unique critical point is given~by
$$
z_{ij}=\frac{a_{++}a_{ij}}{a_{i+}a_{+j}}, \quad \text{where} \quad  a_{i+}=\sum_j a_{ij},\quad
a_{+j}=\sum_i a_{ij},\quad
a_{++}=\sum_{i,j}a_{ij}.
$$
Indeed, the resulting matrix $Z - {\bf 1}$ has rank one because the vector
$(a_{1+},a_{2+})$ is in its left kernel. We further check that the
gradient of the objective function, which is the matrix with entries
$
\frac{a_{ij}}{z_{ij}}
=
\frac{a_{i+}a_{+j}}{a_{++}},
$
lies in the linear span of the gradient matrices of the constraints
$$
(z_{1j}-1)(z_{2k}-1)
-
(z_{1k}-1)(z_{2j}-1).
$$
We have thus identified the unique critical point and written it in
Horn form.
\end{example}

\section{Combinatorics of the ML degree} \label{sec7}

In this section, we discuss the number $\mu(G)$ of complex solutions to the
graphical scattering equations (\ref{eq:gradient}) on a graph $G$.
 We conjecture that the following general formula holds.

\begin{conjecture}\label{thm: MLdeg_admissible}
Let $G$ be a copious and topologically copious graph with $n$
vertices. If the corresponding very affine variety $\mathcal M_G$ is smooth, then its ML degree is given by
\begin{equation}\label{eq:muformula}
\mu(G)
\,\,=\,\,
\biggl|
\sum_{s=3}^n
(-1)^{s-3}(s-3)!\,\pi_G(s)
\biggr|,
\end{equation}
where $\pi_G(s)$ is the number of partitions of $[n]$ into $s$
independent subsets of $G$.
\end{conjecture}
Our computations support Conjecture~\ref{thm: MLdeg_admissible} for all graphs  up to $8$
vertices for which $\mathcal M_G$ is smooth. A previous version of this paper omitted
the smoothness hypothesis. We are grateful to Bernhard Reinke for
pointing out a counterexample. In fact, our database~\cite{BBBBV}
contains $14$ counterexamples to the previous formulation, and the
corresponding very affine varieties $\mathcal M_G$ are singular. For instance, the formula~\eqref{eq:muformula}  fails in general for complete bipartite graphs
$K_{s,t}$ whose very affine variety ${\cal M}_{K_{s,t}}$ admits a determinantal description~by~Theorem~\ref{thm:bipartite}.
\begin{example}[$n=6$]
The $15$ graphs $G$ in Proposition \ref{prop:14} are presented by their
non-edges. The number $\pi_G(s)$ counts
 set partitions of $\{1,2,3,4,5,6\}$ into $s$ blocks,
 where every pair in a block is a non-edge.
 The complete graph $K_6$ has $\pi_G(s) = 0$ for $s \leq 5$ and
 $\pi_G(6) = 1$. The  octahedron graph, 
 with non-edges $12$, $34$ and $56$, has
 $ \,\pi_G(3) = \pi_G(6) = 1$ and $ \pi_G(4) = \pi_G(5) = 3$,
 and the sum in (\ref{eq:muformula}) gives
 $\mu(G) = | 1 \cdot 1 - 1 \cdot 3 + 2  \cdot 3 - 6 \cdot 1| = 2$.
\end{example}

\begin{example}[Two non-edges]
Let $G$ be the complete graph $K_n$ with two edges removed.
If the two non-edges share a vertex, then (\ref{eq:muformula}) gives the formula
$\mu(G) \, =\, (n-3)! - 2(n-4)!$. But,
if the two non-edges of $G$ are disjoint, then
$\mu(G) \,=\, (n-3)! \,-\, 2(n-4)! \,+\, (n-5)!$.
\end{example}

\begin{example}[Removing a triangle]
Let $G$ be the graph on $[n]$ with non-edges $12,13,23$.
Conjecture \ref{thm: MLdeg_admissible} states that
$\mu(G) = (n-3)! - 3(n-4)!+(n-5)!$. 
See Example~\ref{ex:comptriang} for $n=7$.
\end{example}

  \begin{remark}
      Any
     list of $r$      (generalized) cross-ratios, arising from
     lattice vectors $u_1,\ldots,u_r$ in      $\mathcal{K}_{K_n}$
     defines a regular map (\ref{eq:umap2}) from
     $\mathcal{M}_{0,n}$ to $(\CC^*)^r$.
What can we say about this map?
What are its fibers and its image?
The closure of the image is a very affine variety.
We can ask for its ML degree, its  prime ideal, and whether it belongs to a modular
compactification.
If $r=n-3$ and $\tilde u_1, \ldots, \tilde u_r$ are algebraically independent cross-ratios,
then the map is dominant, and  it is $\beta$-to-$1$ where $\beta$ is the cross-ratio degree
studied by Silversmith~\cite{Silv}.
      \end{remark}

\begin{remark}
   For graphs with a universal vertex as in Proposition \ref{prop:gauge}, the 
    partial compactification $\mathcal{M}_G$ coincides with the open part $\mathcal{M}_{0,\Gamma}$ of the moduli space of graphically stable curves 
    studied by Fry \cite{Fry2022}.
    Fry's map extends our embedding for graphs with a universal vertex. 
    See the proof of Proposition \ref{prop:gauge}. Deletion-contraction properties of
these graphically stable spaces are studied in~\cite{KenyonLam,RS, Yang}. For arbitrary graphs many open problems remain.
\end{remark}

The formula in  (\ref{eq:muformula}) is more general than the
known formula for graphs with a universal vertex. When $G$ has a universal vertex, its ML degree is the signed Euler characteristic of the complement of a graphic hyperplane arrangement,
written as in (\ref{eq:gaugefixing}).
By Varchenko's Theorem,
 the number of critical points equals the number of bounded regions in this complement.
By Zaslavsky's Theorem, this number can be determined from the characteristic polynomial of
the underlying matroid. For graphic arrangements, this is the chromatic polynomial.
The {\em chromatic polynomial} $p_H(t)$ of a graph $H$ counts the number 
of colorings of $H$ with $t$ colors. In our situation,
$p_H(t)$ is always divisible by $t-1$, and we write
$\overline{p}_H(t) =  p_H(t)/(t-1)$.
We use this polynomial for the
graph $H = G \backslash n$ that is obtained by deleting the vertex $n$.

\begin{theorem}[Varchenko-Zaslavsky] \label{cor:gauge}
Let $G$ be a copious graph on $[n]$ for which $n$ is a universal vertex, as
in Proposition~\ref{prop:gauge}.
Then the ML degree equals 
\begin{equation}
\mu(G) \,=\, |\,\overline{p}_{G\backslash n} (1)| \,=\,  \biggl| \,\sum_{s=3}^n
(-1)^{\, s - 3} \,\big( s - 3\big)! \,\pi_G(s) \,\biggr|,
\end{equation}
where the righthand side of the second equality is our formula in (\ref{eq:muformula}).
\end{theorem}
\begin{proof}
By Proposition~\ref{prop:gauge}, $G$ is also topologically copious. Hence, we can apply Theorem~\ref{thm:pacom} and identify $\mu(G)$ with the signed Euler characteristic of the very affine variety~$\mathcal{M}_G$.  
In the proof of Proposition~\ref{prop:gauge} we realized $\mathcal{M}_G$ as the complement of the
arrangement $\mathcal{A}_G$ in $\CC^{n-3}$ consisting of the hyperplanes $z_i, z_j - 1, z_j - z_i$,
where $1i$, $2j$, and $ij$ are edges of the graph $G\backslash n$.  
By Remark \ref{rem:deconing}, this arrangement is  the deconing of the graphical arrangement $\mathcal{H}_{G\backslash n}$.  
By Varchenko’s Theorem, the number of critical points  equals the number of bounded regions in the corresponding real arrangement complement. The latter number, up to sign, is given by evaluating the characteristic polynomial  $\chi(\mathcal{A}_G, t)$ 
of the hyperplane arrangement at $t=1$.  

We now apply the relationship, found in \cite[Proposition 3.1.2]{OrlikTerao2001},
 between the characteristic polynomial of the deconed arrangement and that of the original graphical arrangement:
\[ \chi(\mathcal{A}_G, t)
\,\,=\,\, \frac{\chi(\mathcal{H}_{G\backslash n}, t)}{t-1}
\,\,=\,\, \frac{p_{G\backslash n}(t)}{t-1}
\,\,= \,\,\overline{p}_{G\backslash n}(t).
\]
Here we used the classical fact that the characteristic polynomial of a graphical arrangement $\mathcal{H}_G$ coincides with the chromatic polynomial $p_G(t)$ of the graph~$G$. From this we conclude 
$\,
\mu(G)
= \bigl|\,\chi(\mathcal{M}_G)\,\bigr|
= \bigl|\,\chi(\CC^{n-3}\backslash \mathcal{A}_G)\,\bigr|
= \bigl|\,\chi(\mathcal{A}_G, 1)\,\bigr|
= \bigl|\,\overline{p}_{G\backslash n}(1)\,\bigr|$.
This is our first equality.

Since $n$ is universal, notice that $\pi_G(s)=\pi_{G\backslash n}(s-1)$, and thus $p_{G\backslash n}(t)=\sum_{r=1}^{n-1} (t)_{r}\pi_{G\backslash n}(r),$ where $(t)_r=t(t-1)\cdots (t-r+1).$ We then compute the chromatic polynomial as follows:
\begin{align*}
    \overline{p}_{G\backslash n}(1)\,\,=\,\,\sum_{r=1}^{n-1} \,(1-2)(1-3)\cdots (1-r+1)\pi_{G\backslash n}(r)\,\,=\,\,\sum_{r=1}^{n-1}(-1)^r (r-2)!\,\pi_{G\backslash n}(r).
\end{align*}
Setting $s=r+1$ and taking the absolute values  produces our second equality.
\end{proof}

\begin{example}[$n=6$]
Among the $15$ graphs in Proposition \ref{prop:14},  only ten
satisfy the hypothesis of Theorem \ref{cor:gauge}.
If $G$ is the graph with one non-edge $12$,
then $G \backslash 6$ has nine edges, and its chromatic polynomial gives
$\overline{p}_{G\backslash 6}(t) = t(t-2)(t-3)^2$. If $G$ has three non-edges $12,34,45$
then $\overline{p}_{G\backslash 6}(t) = t(t-2)(t^2-4t+5)$.
The evaluations at $t=1$ are $-4$ and $-2$.
Up to sign, these are the numbers seen in the column $\mu(G)$ 
for the table in Proposition \ref{prop:14}.
\end{example}
\section{Graphs on Few Vertices and Hypertrees} \label{sec8}

The previous sections introduced four different 
ways for a graph $G$ to be copious.
Theorem~\ref{conj:main} tells us that three are equivalent. 
We continue to call $G$ {\em copious} if it satisfies the equivalent conditions in
Definitions \ref{def:geomcop}, \ref{def:matrcop} and \ref{def:algcop}.
This section is dedicated to two topics. The first is the findings of our comprehensive computational study of graphs with up to nine vertices.
 The blueprint is given by the
lists of copious graphs for $n=5$ in Example \ref{ex:ideals5}  and for $n=6$ in
 Proposition \ref{prop:14}.
We computed the analogous lists in the much harder cases $n=7,8$, with partial results also for $n=9$. Second, we investigate topological copiousness for hypertrees.

The notion of a hypertree was introduced by Castravet and Tevelev in
 \cite{CT}, and its connection to scattering equations was featured in \cite[Section 5]{EPS}. 
 We note that
 Questions 5.5 and 5.7  in \cite{EPS}
  were our  original motivation for launching this project.
  But,  our investigations then took a different route, and these
   two questions remain largely open.

From the entry A000088 in the Online Encyclopedia of Integer Sequences (OEIS), we know that
the number of graphs 
for $n=4,5,6,7,8,9$ equals
$ 11, 34, 156, 1044, 12346, 274668$.
The data for the following theorem appears on our
supplementary materials website \cite{BBBBV}.

    \begin{theorem}\label{thm:copious-graphs}
    The number of copious graphs for $n=4,5,6,7,8$ equals $1,3,15, 129, 2328$.
        For each of these graphs, the map~\eqref{eq:umap2} is generically
injective, which is necessary for being topologically copious. For $n=9$, there are at least $26324$ copious graphs.    
    \end{theorem}

    \begin{proof}[Proof and Discussion]
        We obtain the lists of all graphs on $n=6,7,8,9$ vertices up to isomorphism from the graph library
in {\tt SageMath}.  To identify the copious graphs, we first restrict to graphs that are compatible
with momentum conservation. This means that the
even cycle matroid $M(\mathcal{K}_G)$ has no co-loop.
We then compute the scattering matroid $\mathcal{S}(G)$
and check that it has no co-loops. Such a co-loop would
 prevent the map $\mathcal{V}_G \rightarrow \mathcal K_G$ from being dominant. 
 At this stage we also check that the rank of $\mathcal{S}(G)$ satisfies inequality (\ref{eq:ScatMatrank}).
 
 We wish to decide, for each
 remaining graph $G$, whether $G$ is matroidally copious.
For this we need to check the equivalent conditions (\ref{eq:fromxtou2})–(\ref{eq:fromxtou1}).
We do this by a probabilistic method. Namely, we pick a 
reasonably large sample
of rational points $x \in\mathbb{P}^{n-1}$. 
If $\kappa(G)$ equals the dimension on the left hand side of (\ref{eq:fromxtou1}) for our
sample, then we are done. Otherwise, we have strong 
evidence that $G$ is not copious.
If so, it is best to try again with a larger~sample.
 
To gain certainty, we turn to the algebraic definition of copious.
  We aim to compute the multidegree (\ref{eq:multidegree}) for each graph $G$.
  We exclude $G$ whenever
  the ratio in (\ref{eq:muG}) is not a positive integer.
    For $n=6$, we obtain all multidegrees symbolically with {\tt Macaulay2}~\cite{Macaulay2Book} within a few hours. For $n=7$, a significant 
    fraction of the graphs require days. We therefore perform this task with numerical methods.
    Namely, we use the {\tt Julia} package {\tt HomotopyContinuation.jl}~\cite{HomCont}.
    This allows us to compute the multidegrees for $n=8$~as~well. 
    
        The steps described in the first paragraph above yield $140$ graphs with
                 co-loop-free scattering matroid for $n=7$. Among these, $11$ have a multidegree 
                 whose first  coefficient is not divisible by $5$. This leaves us with $129$ graphs.
        To certify that these $129$ are matroidally copious, we checked (\ref{eq:candidate}) as well as (\ref{eq:fromxtou1}).
        This was done in {\tt Macaulay2}~\cite{Macaulay2Book} with $21$ points $x$ sampled from~$\mathbb{P}^{6}$.
        The same methodology also worked well for $n=8$ and $9$, giving a lower bound for the number of matroidally copious graphs. For $n=4,5,6,7,8$, the matroidal computations certify that the respective lists contain $1,3,15,129,2328$ copious graphs. In addition, our computations show that the map~\eqref{eq:umap2} is
generically injective for each graph~in~these~lists. 
        
       To verify generic injectivity of the map~\eqref{eq:umap2}, we again use
\texttt{Macaulay2}~\cite{Macaulay2Book}. 
         For a random point $x\in \mathcal{M}_{0,n},$ we check that the preimage 
        of $(\tilde{u}_1(x),\tilde{u}_2(x),\dots,\tilde{u}_r(x))$ 
        under the map (\ref{eq:umap2}) is a single {\rm PGL(2)} orbit.
        This is done by computing  an elimination ideal, but after
                  fixing $x_1=1,\; x_2=2,\; x_3=3$. The resulting ideal for
                                      the preimage must then be a maximal ideal.

         Our method for identifying all matroidally copious graphs runs well also for larger $n$. But, it has the
         disadvantage that it is probabilistic.
                    In practice, we can sample only finitely many $x \in\mathbb{P}^{n-1}$ when certifying 
                    the linear algebra conditions (\ref{eq:fromxtou2}) and (\ref{eq:fromxtou1}).
                      By contrast, we found that computing multidegrees is infeasible
                      for $n=9$ with the numerical method.
      \end{proof}

We next present a more detailed view of the data in Theorem \ref{thm:copious-graphs}.
Table \ref{tab:78} counts copious graphs $G$ for fixed
$|G|$ and fixed $ \mu(G) $. Each row and column corresponds to a
range of values.

\vspace{-0.77cm}

\begin{center}
\begin{table}[h]
$$ \vspace{-0.4cm} \!\!\!\!\! \!\begin{small} \begin{matrix}   
\qquad |G| \backslash  \mu(G)\!\! \!\!& 1 & 2 & 3 & 4 & \!5,\!6\! &\!\! 7 .. 9 \!\!& \!\!10 ..24\! \! \\
\!\! 12 & \! 10 & 2 & 0 & 0 & 0 & 0 & 0 \\
\!\! 13 & \! 14 & 8 & 1 & 0 & 0 & 0 & 0 \\
\!\! 14 & \! 11 & 13 & 5 & 2 & 0 & 0 & 0 \\
\!\! 15 & 5 & 9 & 5 & 6 & 3 & 0 & 0 \\ 
\!\! 16 & 0 & 4 & 2 & 5 & 6 & 0 & 0 \\
\!\! 17 & 0 & 0 & 0 & 3 & 2 & 4 & 0 \\
\!\! 18..21 & 0 & 0 & 0 & 0 & 1 & 2 & 6 
\end{matrix} \end{small}
\qquad \,\,\,\,
\begin{footnotesize}
 \begin{matrix}   
\quad |G| \backslash  \mu(G)\!\! \!\!&
1 & 2 & 3,\! 4 & 5,\! 6 & \!\! 7..10 \!\! & \!\! 11..20 \!\! & \!\! 21..120 \!\! \\
13 & 0  & 1  & 0  & 0  & 0  & 0  & 0 \\
 14 & 52  & 29  & 3  & 0  & 0  & 0  & 0 \\
 15 & \! 109  & 95  & 39  & 1  & 0  & 0  & 0 \\
16 & \! 104  & \!156  & \! 126  & 20  & 2  & 0  & 0 \\
17 &  52  &\! 133  & \! 185  & 89  & 15  & 0  & 0 \\
18 &  15  & 65  & \!  142  & 103  & 96  & 5  & 0 \\
19 &  0  & 16  & 63  & 63  & 114  & 50  & 0 \\
20  , \! 21 &  0  & 0  & 17  & 24  & 75  & 154  & 20 \\
22..28 &  0  & 0  & 0  & 0  & 0  & 23  & 72
\end{matrix}
\vspace{-0.2cm}
\end{footnotesize}
$$
\caption{Copious graphs by number of edges and ML degree for $n\!=\!7$ (left) and $n\!=\!8$ (right).
\label{tab:78} }
\end{table}
\end{center}

\vspace{-0.67in}

The first row for $n=8$ shows a unique graph $G$ with $|G| = 13$.
Interestingly, it has $\mu(G) = 2$.
This $G$ is the complete bipartite graph $K_{4,4}$ with three
disjoint edges removed.

In the first column we see copious graphs $G$ with ML degree $\mu(G) = 1$.
If the edge set of $G$ is inclusion-maximal then this is a case of
minimal kinematics \cite{EPS}.
The study in \cite{EPS} was restricted to graphs with a universal vertex, with
coordinates as in (\ref{eq:gaugefixing}).
Removing that universal vertex, the relevant graphs $H$ on $[n-1]$
are precisely the $2$-trees \cite[Theorem 1.4]{EPS}.

According to A054581 in the OEIS, there are twelve  $2$-trees
for $n=8$. Each of them corresponds to a graph $G$ with 
$3n-6 = 18$ edges. In  Table \ref{tab:78} we see $15$ graphs with
$|G| \!=\! 18$ and $\mu(G) \!=\! 1$. Here is one of the three minimal kinematics
graphs not coming from a~$2$-tree.

    \begin{example}
        Let $n=8$ and let $G$ be the graph with non-edges $12$, $13$, $14$, $18$, $23$, $24$, $27$, $34$, $36$, $45$. Since $G$ has no universal vertex, it does not arise from a 2-tree and its multidegree~is 
        \[ 6\sigma^8\tau^5 + 31\sigma^9\tau^4 + 52\sigma^{10} \tau^3+ 36\sigma^{11}\tau^2 + 10\sigma^{12} \tau + \sigma^{13}.\]
        Dividing the leading coefficient by $n-2=6$, 
                equation (\ref{eq:muG}) tells us that $\mu(G)=1$.
        \end{example}

We state  a conjecture which is consistent with our data and with the results in \cite{EPS}.

\begin{conjecture}
Every minimal kinematics graph on $[n] = \{1,2,\dots,n\}$ has $3n-6$ edges.
\end{conjecture}
Next comes a graph which is counted in the entry $2$ in the last row on the left in Table \ref{tab:78}.

    \begin{example} \label{ex:comptriang}
        Let $n=7$ and consider the graph $G$ with non-edges $12,13,23$.
        This graph is copious with ${\rm rank}(\mathcal{S}(G)) = 11 =7+(7-3)<18=|G|$.
         We describe the embedding (\ref{eq:umap2}) of        
         $\mathcal{M}_{0,7}$ into $(\mathbb{C}^*)^{11}$. 
        We write $[ij|k\ell]=\frac{[ik][j\ell]}{[i\ell][jk]}=\frac{(x_i-x_k)(x_j-x_{\ell})}{(x_i-x_{\ell})(x_j-x_k)}$ for the usual cross-ratio.
                For one choice of basis of the kinematic lattice $\mathcal{K}_G \cap
         \ZZ^{|G|} = \mathbb{Z}^{11}$, the embedding  is given by
        $$ \begin{small}
         \left([12|45],[12|46],[12|47],[13|45],[13|46],[13|47],[14|56],[14|57],[15|46],[15|47],
         \frac{[14][15][56]}{[16][17][45]}\right)\!.
         \end{small} $$
        The multidegree of $\overline{\mathcal{V}}_G$ equals 
        $40\sigma^7\tau^4 + 74\sigma^8 \tau^3+ 45\sigma^9\tau^2 + 11\sigma^{10}\tau + \sigma^{11}$, and hence $\mu(G)=8$.
    \end{example}

    Our supplementary materials website \cite{BBBBV} contains all
  graphs for $n=7,8,9$, together with our \texttt{Julia} and \texttt{Macaulay2} code, as well as the datasets of ML degrees and multidegrees of scattering ideals. This is posted at the
open-access  repository {\tt Zenodo} for research data:
    \begin{center}
    \url{https://zenodo.org/records/17546339}.
    \end{center}

Our second topic is hypertrees \cite[Section 5]{EPS}. We identify each hypertree $T$ with a graph $G$. Its edges are the pairs
 in the triples of $T$.
The smallest such hypertree graph $G$ is the octahedron \cite[eqn~(26)]{EPS}, which belongs to a special family of hypertrees called bipyramids.  
The {\em bipyramid graph} $G$ has $n$ vertices
and $3n-6$ edges. Namely, we start with the $(n-2)$-cycle labelled by
$3,4,5,\ldots,n,3$, and we connect both $1$ and $2$ with each vertex in the cycle. 
The following result can be viewed as a step towards answering \cite[Question 5.5]{EPS}.

\begin{proposition}\label{prop:bipyramidTopologicallyCopious}
Every bipyramid is topologically copious.
\end{proposition}
\begin{proof}
Let $G$ be the bipyramid on $[n]$ and write
 $\{e_{ij}\}_{ij\in G}$ for the standard basis  of $\mathbb{Z}^{3n-6}$.
  Then $\{e_{13}+e_{2i}-e_{1i}-e_{23}\}_{i=4}^n$ is a linearly independent set in 
  $\mathcal{K}_G$. 
 These vectors form a primitive set in the kinematic lattice and can therefore be completed to a $\mathbb Z$-basis
  $\{e_{13}+e_{2i}-e_{1i}-e_{23}\}_{i=4}^n\cup \{{u}_{n-2},\dots,{u}_r\}$ of $\mathcal{K}_G\cap\ZZ^{3n-6}$.
Let $\varphi:\mathcal{M}_{0,n}\to(\mathbb{C}^\ast)^r$ be the corresponding map (\ref{eq:umap2}).
The generalized cross-ratios given by the basis elements $e_{13}+e_{2i}-e_{1i}-e_{23}$  have~the~form 
\[\frac{(x_1-x_3)(x_2-x_i)}{(x_1-x_i)(x_2-x_3)} \,\,=\,\,[12|3i] \,=\,[21|i3]. \]
 After fixing $x_1=\infty,x_2=0,x_3=1$, we find that $[21|i3]=x_i$.
 Hence the affine coordinates $x_4,\ldots,x_n$ on
$\mathcal M_{0,n}$ occur among the coordinates of $\varphi$.
Projection onto $x_4,\ldots,x_n$ gives a regular inverse map
of $\varphi$. Therefore, $\varphi$ is an embedding, and $G$ is
topologically~copious.
\end{proof}


Moving beyond bipyramids, we believe this result holds for all irreducible hypertrees. 

\begin{conjecture}
    All irreducible hypertrees are topologically copious.
\end{conjecture}

To this end, we certified using {\tt  Macaulay2}~\cite{Macaulay2Book} that the map \eqref{eq:umap2} is generically injective for all irreducible hypertrees with $n\leq 10$ vertices in the Opie-Scheidwasser database 
    \begin{center}
    \url{https://websites.umass.edu/tevelev/hypertree_database/}.
    \end{center} 
    Concretely, for a random
point $x\in \mathcal{M}_{0,n}$, we check that the preimage of~$(\tilde{u}_1(x), \ldots,\tilde{u}_r(x))$ under the map \eqref{eq:umap2}
is a single $\text{PGL}(2)$ orbit. This is done by fixing 
$x_1=1,\; x_2=2,\; x_3=3$ and computing an elimination ideal. The resulting ideal for the preimage must then be a maximal ideal. This verifies generic injectivity, which is weaker than the embedding property. 

We conclude this section with a conjecture for the ML degree of the bipyramid graph $G$ on $n$ vertices. For $6\leq n\leq 15,$ we confirmed that $\mu(G)$ is $n-4,$ and suspect this for all $n$.

\begin{conjecture}\label{conj:bipyramid}
    If $G$ is a bipyramid on $n$ vertices, then $\mu(G)=n-4.$
\end{conjecture}

Happily, this coincides with the 
alternating sum in (\ref{eq:muformula}).

\begin{proposition}\label{thm:bipyramid}
If $G$ is a bipyramid on $n$ vertices, then
(\ref{eq:muformula}) is equal to $n-4$.
\end{proposition}

\begin{proof}
Any partition of $ [n]$ into $s$ independent subsets either has $1$ and $2$ together in a part of size two or separately in two parts of size one. The number of such partitions is thus counted by partitions of the 
vertices of the $(n-2)$-cycle into $s-1$ and $s-2$ independent subsets. 
The chromatic polynomial of the $m$-cycle is $(t-1)^m+(-1)^m(t-1)$, and the number of partitions of a graph $H$ into $k$ independent sets equals $p_H(k)/k!$. This implies
\begin{align*}
        \pi_G(s)&\,\,=\,\,\, \frac{1}{(s-1)!}\sum_{i=0}^{s-1}(-1)^{s-1-i}\binom{s-1}{i}((i-1)^{n-2}+(-1)^{n-2}(i-1))\\&\quad +\,\,\frac{1}{(s-2)!}\sum_{j=0}^{s-2}(-1)^{s-j}\binom{s-2}{j}((j-1)^{n-2}+(-1)^{n-2}(j-1)).
\end{align*}
The absolute value of the nested sum that arises by
substituting $\pi_G(s)$ into (\ref{eq:muformula}) equals $n-4$. A computer-aided proof was found by Manuel Kauers.
We present this on our supplementary  website \cite{BBBBV}.
\end{proof}

\section{Discriminants} \label{sec9}

Our last section is devoted to the {\em logarithmic discriminant},
which vanishes when the scattering potential $L_G$ has fewer than $\mu(G)$ critical points.
This was studied by Kayser,  Kretschmer and Telen \cite{KKT}
 for the complete graph $ K_n$ and it is here extended to any graph~$G$.
 
Let $G$ be a copious graph with $\mu(G) \geq 2$.
The  {\em logarithmic discriminant} $\Delta_{G}= \Delta_G(s)$ is an irreducible
homogeneous polynomial of positive degree $\delta(G)$ in the
unknowns $s_{ij}$.
The hypersurface of $\Delta_G$ 
is the branch locus of 
the $\mu(G)$-to-$1$ map $(s,x) \mapsto s$ from the
scattering correspondence $\mathcal{V}_G$ onto the kinematic space $\mathcal{K}_G$.
Note that the  logarithmic discriminant $\Delta_G$  is not unique as a polynomial
in the $s_{ij}$. It
is unique only
modulo the linear relations~(\ref{eq:momcon}).

We begin with $n=6$.
The complete graph has $\mu(K_6) = 6$, and the degree of its
discriminant is $\delta(K_6) = 30$. This was found numerically
by Kayser, Kretschmer and Telen~\cite{KKT}. For
details see Conjecture~1 and Example 8.5 in \cite{KKT}.
The value of $\delta(G)$ for other graphs $G$ is shown
in our Proposition~\ref{prop:14}.
We now focus on the  octahedron
(Examples \ref{ex:octah} and \ref{ex:octah3}).

The $12$ Mandelstam parameters $s_{ij}$
must satisfy momentum conservation at each of the
six vertices, and there are many different choices
of bases modulo these relations. There are four distinguished basis choices
which respect the symmetry among the six vertices of $G$.
Namely, we fix one of the four pairs of antipodal
facets of the octahedron. Then the six edges not
on those facets form a $6$-cycle, and these
six $s_{ij}$ form a basis.  Here is one of them:
\begin{equation}
\label{eq:secondbasis} \begin{footnotesize} \begin{pmatrix}
s_{14}  \\ s_{16} \\ s_{46}  \\ s_{23} \\ s_{25} \\ s_{35} \\
\end{pmatrix}
\quad = \quad
\frac{1}{2} 
\begin{pmatrix}
          -1  &  -1  &  -1  &  -1  & \phantom{-} 1  & \phantom{-} 1 \,\, \\
          -1  &  -1  &  \phantom{-} 1  &  \phantom{-} 1   &  -1  &  -1 \,\, \\
          \phantom{-} 1  &  \phantom{-} 1  &  -1  &  -1  &  -1  &  -1 \,\, \\
          -1  &  \phantom{-} 1  &  \phantom{-} 1  &  -1  &  -1  &  -1 \, \, \\
          \phantom{-} 1  &  -1  &  -1  &  -1  &  -1  &  \phantom{-} 1 \,\, \\
          -1  &  -1  &  -1  &  \phantom{-} 1  &  \phantom{-} 1  &  -1 \,\,
\end{pmatrix}          \!
\begin{pmatrix}
s_{13} \\ s_{15}  \\ s_{45}  \\ s_{24} \\ s_{26} \\ s_{36}
\end{pmatrix}.
\end{footnotesize}
\end{equation}
These six linear equations are equivalent to (\ref{eq:momcon}).
In the left column vector we see the antipodal facets
$146$ and $235$, and the right column vector is the
$6$-cycle $1-5-4-2-6-3-1$.

To compute $\Delta_G$, we 
start with the equations $\partial L_G /\partial x_i = 0$.
Substituting (\ref{eq:secondbasis}),
we obtain a system of bihomogeneous
equations in $6$ unknowns $s_{ij}$
and $6$ unknowns $x_i$. 
Over a general point in $\mathcal{K}_G$,
these equations define the union of two
surfaces in $\PP^5$ as in (\ref{eq:I4}).
By specializing three of the coordinates,
say $x_4 = 1, x_5 = 2, x_6 = 3$,
we obtain two points.
We eliminate $x_2,x_3$ from the system and
obtain one quadratic equation in $x_1$
whose coefficients are expressions in the
six parameters $s_{ij}$.
The discriminant of this quadric in $x_1$ has
the logarithmic discriminant $\Delta_G$ as
its main factor. Extraneous factors are
removed by changing the variable
order in the elimination. Namely, $\Delta_G$ is
the gcd of the discriminants of the three quadrics in $x_1,x_2,x_3$.

We find that the discriminant $\Delta_G$ is a homogeneous polynomial
of degree six on $\mathcal{K}_G$:
$$ \Delta_G\, = \,
s_{13}^6+2 s_{13}^5 s_{15}+2 s_{13}^5 s_{24}-2 s_{13}^5 s_{26}+2 s_{13}^5 s_{36}
-2 s_{13}^5s_{45}-s_{13}^4 s_{15}^2+6 s_{13}^4 s_{15} s_{24}-6 s_{13}^4 s_{15}s_{26}
+\,\cdots. $$
This is a sum of $462$ monomials in the basis (\ref{eq:secondbasis}).
This discriminant has remarkable properties. If we set one of the variables
$s_{ij}$ to zero, then $\Delta_G(s)$ is the square of a cubic with $35$ terms.
If we set two variables $s_{ij}$ to zero, then it becomes the
product of three squares:
$$ \Delta_G|_{s_{13} = s_{45} =0} \,\,= \,\,
(s_{15} - s_{24} + s_{26} + s_{36})^2 \cdot  (s_{15} - s_{24} - s_{26} - s_{36})^2  \cdot (s_{15} + s_{24} + s_{26} - s_{36})^2. $$
The discriminant $\Delta_G$ takes on both positive and negative values. For instance, 
$\,\Delta_G(s) = -71 \,$ for $\,s_{13} = s_{15} =  s_{24} =  s_{26} =  s_{45} = 1, s_{36} = -2$.
Here, both critical points of $L_G$ are non-real.

We found experimentally that
the logarithmic discriminant $\Delta_G(s)$ of
the octahedron $G$ is {\em copositive} when written
 in the $6$-cycle basis (\ref{eq:secondbasis}) for $\mathcal{K}_G$, i.e.,
$\Delta_G(s) > 0$ for all $s \in \RR_{>0}^6$.
What makes this observation remarkable is that there is no gauge fixing
for the octahedron. For all graphs $G$ with a universal vertex,
copositivity follows immediately from \cite[Corollary 4.2]{KKT}.
For instance, for $K_5$, this is the AM-GM inequality in \cite[Example 1.1, page 3]{KKT}.

\begin{problem}
Study the copositive geometry \cite{Mate} of the logarithmic discriminant
$\Delta_G(s)$ for all graphs $G$ that do not satisfy the
 hypothesis in Proposition \ref{prop:gauge} and
 Theorem \ref{cor:gauge}.
\end{problem}

 We now discuss  numerical computations of the
discriminant degrees $\delta(G)$ for $n=6,7,8$.
 From \cite{KKT} we already know that
 $\delta(K_6)=30,\; \delta(K_7)=208,\; \delta(K_8)=1540$.
 For simplicity, we restrict our computations to those graphs that have
 a universal vertex, as in Proposition~\ref{prop:gauge}.
 
  \begin{proposition} \label{prop:deltadata}
      For the copious graphs $G$ on $[n]$ with a universal vertex,  we found  $\delta(G)\in \{4,10,16,30\}$ for $n=6$, 
             $\delta(G)\in \{4,8,10,16,18,24,30,32,40,44,48,66,76,$ $ 84,104,144,208\}$ for $n=7$, and
       $4\leq \delta(G)\leq 1540$ for $n=8$, where $\delta(G)$ was always even.
  \end{proposition}
  
  \begin{proof}
  We use the {\tt Julia} code due to
      Kayser,  Kretschmer and Telen \cite{KKT}
      for logarithmic discriminants. This is found at
\url{https://mathrepo.mis.mpg.de/LogDiscHyperplanes/}.
The computation for the complete graph $G = K_n$ is described in \cite[Section 8]{KKT}.
It rests on the matrix $\mathcal{A}$ which encodes the braid arrangement, in the
affine version of Remark \ref{rem:deconing}.

For any graph $G$ with a universal vertex $n$, we remove      
        the  columns of the matrix $\mathcal{A}$ that are indexed by the non-edges of $G$.
                 In some cases, we first performed a permutation of $[n]$ to ensure that the $n$
         pairs $12,\,1n,2n,\ldots,(n-1)n$ are edges of $G$.
         This is necessary because~\cite{KKT} fixes the gauge $x_1=0,\; x_2=1,\; x_n=\infty$.
         Let ${\tt s}$ denote the resulting set of non-edges of $G$.
         This is a subset of $\binom{[n-1]}{2}$. We now run the following Julia code within the package of \cite{KKT}:
      \begin{verbatim}
      L = get_L_M0m(n)
      L' = L[:, setdiff(1:size(L,2),s)]
      get_deg_DeltaLog(L')
      \end{verbatim}
      \vspace{-0.4cm}
      The output is the  discriminant degree $\delta(G)$ of the copious graph $G$. The computation
            occasionally led to a numerical error, in which case
            we had to redo it more carefully.
  \end{proof}


\bigskip 

\noindent {\bf Acknowledgement}:
We thank Manuel Kauers for the proof of Proposition \ref{thm:bipyramid} and Benjamin Hollering for invaluable input on computations and for  the list of graphs on $n=9$ vertices. We are grateful to Veronica Calvo Cortes, Jake Levinson, and Rob Silversmith for insightful discussions and to Bernhard Reinke for  a counterexample to a conjecture that appeared in a previous version of this paper. We thank the anonymous referee for the valuable comments and suggestions.
BF and BZ 
were both supported by Graduate Research Fellowships in mathematics from the
US National Science Foundation.
VB and BS acknowledge support
from the European Research Council \begin{small} (UNIVERSE PLUS, 101118787). \end{small}
$\!\!$ \begin{scriptsize}Views~and~opinions expressed
are however those of the authors only and do not necessarily reflect those of the European Union or the 
European
Research Council Executive Agency. Neither the European Union nor the granting authority
can be held responsible for them.
\end{scriptsize}

\bigskip 

\noindent
\footnotesize {\bf Authors' addresses:}
\smallskip

\noindent Barbara Betti, OvGU Magdeburg, Germany \hfill {\tt  barbara.betti@ovgu.de}

\noindent Viktoriia Borovik, MPI-MiS Leipzig, Germany \hfill {\tt  borovik@mis.mpg.de}

\noindent Bella Finkel, U Wisconsin, Madison, USA \hfill {\tt blfinkel@wisc.edu}

\noindent Bernd Sturmfels, MPI-MiS Leipzig, Germany \hfill {\tt bernd@mis.mpg.de}

\noindent Bailee Zacovic, U Michigan, Ann Arbor, USA \hfill {\tt bzacovic@umich.edu}

\end{document}